\documentclass[12pt]{amsart} 
\usepackage{amsmath,amssymb}
\usepackage{latexsym}
\usepackage{amssymb} 
\usepackage{latexsym}
\usepackage{delarray}

\numberwithin{equation}{section}
\theoremstyle{plain}
\newtheorem{thm}{Theorem}[section]
\newtheorem{lem}[thm]{Lemma}
\newtheorem{prop}[thm]{Proposition}

\newtheorem{rem}[thm]{Remark}

\newcommand\R{{\mathbb R}}

\newcommand\N{{\mathbb N}}

\newcommand\Rn{{{\mathbb R}^n}}
\newcommand\Snm{{{\mathbb S}^{n-1}}}

\def\p#1{{\left({#1}\right)}}
\def\jp#1{{\left\langle{#1}\right\rangle}}
\def\abs#1{{\left|{#1}\right|}}
\def\br#1{{\left[{#1}\right]}}
\def\brfig#1{{\left\{{#1}\right\}}}
\newcommand\Ivar{\nu} 

\newcommand\al{\alpha}

\newcommand\ga{\gamma}
\newcommand\de{\delta}

\newcommand\ka{\kappa}
\newcommand\la{\lambda}
\newcommand\va{\varphi}

\newcommand\Si{\Sigma}
\newcommand\om{\omega}

\newcommand\pa{\partial}
\def\Fcal{\mathcal{F}}
\def\supp{\,\textrm{\rm supp}\,}
\def\Re{{\rm Re}}
\def\Im{{\rm Im}}

\title[Dispersive estimates for hyperbolic equations]
{Asymptotic integration and dispersion for hyperbolic
equations}
\author[Tokio Matsuyama and Michael Ruzhansky]
{Tokio Matsuyama${}^\dagger$ and Michael Ruzhansky${}^*$} 
\address{ 
${}^\dagger$Department of Mathematics \endgraf 
Tokai University \endgraf 
Hiratsuka \endgraf 
Kanagawa 259-1292 \endgraf 
Japan \endgraf 
{} \endgraf 
\bigskip
${}^*$Department of Mathematics \endgraf
Imperial College London\endgraf
180 Queen's gate \endgraf 
London SW7 2AZ \endgraf 
United Kingdom} 
\thanks
{ 
\endgraf 
The first author was supported by Grant-in-Aid for 
Scientific Research (C) 
(No. 21540198), Japan Society for the Promotion of Science. 
The second author was supported by the Leverhulme 
Research Fellowship and by the EPSRC grants EP/E062873/1
and EP/G007233/1.
} 
\email{tokio@keyaki.cc.u-tokai.ac.jp \endgraf 
m.ruzhansky@imperial.ac.uk} 
\begin{document}

\maketitle 

\begin{abstract} 
The aim of this paper is to establish time decay properties
and dispersive estimates for strictly 
hyperbolic equations with homogeneous symbols
and with time-dependent coefficients 
whose derivatives belong to $L^1(\mathbb{R})$. 
For this purpose, the method of asymptotic integration
is developed for such equations and representation
formulae for solutions are obtained. These formulae
are analysed further to obtain time decay of $L^p$--$L^q$
norms of propagators for the corresponding Cauchy problems.
It turns out that the decay rates can be expressed in terms
of certain geometric indices of the limiting equation and we
carry out the thorough analysis of this relation. 
This provides a comprehensive view on asymptotic
properties of solutions to time-perturbations of hyperbolic 
equations with constant coefficients.
Moreover, we also
obtain the time decay rate of the $L^p$--$L^q$ 
estimates for equations of these kinds, so
the time well-posedness of the corresponding nonlinear
equations with
additional semilinearity can be treated by standard Strichartz
estimates.
\end{abstract}



\section{Introduction}
\setcounter{equation}{0}
This paper is devoted to several aspects of strictly 
hyperbolic equations of higher orders or of strictly hyperbolic
systems with time-dependent coefficients. In particular,
we will investigate the following topics:
\begin{itemize}
\item representation of solutions of equations of higher order;
\item dispersive estimates for solutions.
\end{itemize}
Equations of orders larger than two appear
often in the analysis of large first order
systems and in the analysis of coupled equations of higher
orders. In the present paper we will restrict our attention
to the investigation of equations with homogeneous symbols
(at the same time making a suitable preparation for the 
further development
of this topic for equations with low order terms).
In fact, we will concentrate on scalar equations of some
order $m\in\N$ keeping in mind that in the case of a system
its dispersion relation (the determinant) will be of such
form, so the information on solution to the Cauchy problem for
the dispersion relation will imply the information on 
solutions to the Cauchy problem of the original system.
On one hand, we will introduce several techniques
allowing to deal with equations of higher orders. 
On the other hand, already for
the second order equations the new method that we propose in this
paper will yield certain improvements and extensions
of known results. In particular,
we will improve the result on the decay
rates in the dispersive estimates
already for the time-dependent wave equation, as well as 
the time decay rate for the 
standard Kirchhoff equation, thus also improving the corresponding
Strichartz estimates. It will also allow the inclusion of mixed
terms in second order equations (a question which is known to
be very delicate if we want to treat problems outside the
perturbation framework).
We will allow time-dependent coefficients and will assume
that their derivatives are in $L^1(\R)$. It is known that
this property is satisfied in many situations, for example
in applications to Kirchhoff equations and systems, etc.
This will also allow us to obtain a comprehensive view on 
time-perturbations 
of equations with constant coefficients, in which case the
assumption of the integrability of derivatives of
coefficients is quite natural.
For the purposes of this paper, we will develop the
{\em asymptotic integration method} for hyperbolic partial 
differential equations with time-dependent coefficients.
While this method is relatively
well-known in the theory of ordinary differential equations
(see e.g. Hartman \cite{Hartman}),
its use in partial differential equations appears to be new.

We note that equations with constant coefficients have
been thoroughly studied by Sugimoto in a series of papers 
\cite{sugi94,sugi96,sugi98} who described
several interesting geometric quantities responsible for the
rate of the time decay of $L^p$--$L^q$ norms of their propagators.
In particular, one has to look at the level
sets of the characteristic roots of the symbol and at the
orders with which tangent lines touch these sets. These
orders become responsible for the time decay rate in
the corresponding dispersion estimates and for
indices of the subsequent Strichartz estimates.
In fact, the appearing indices are related to the oscillation
indices of integral kernels of the propagators, viewed as
oscillatory integrals, and their classification is well 
studied in the singularity theory (e.g. \cite{AVG}).

The case of strictly hyperbolic equations with constant
coefficients with lower order terms has been thoroughly investigated
in \cite{RS2}. In particular, properties of characteristic
roots are crucial in determining exact decay rates and the
complete analysis is quite lengthy and involved. For example,
in the case of equations of dissipative types analysed
in \cite{RS} the decay rate is determined by properties of
characteristics for small frequencies. A general analysis of
this type is necessary for application to large systems,
such as Grad systems in gas dynamics, or to Fokker-Planck
equations, in which case the Galerkin approximation produces
a sequence of scalar equations with orders going to infinity,
see e.g. \cite{R}.
Applications to such problems give a strong additional motivation
to the investigation of equations of higher orders of the
type of those treated in this paper.
 
To become more precise, 
we consider the Cauchy problem for an $m^{\rm th}$ 
order strictly hyperbolic 
equation with time-dependent coefficients, for 
function $u=u(t,x)$: 
\begin{equation}
L(t,D_t,D_x)u 
\equiv D^{m}_t u+\sum_{\underset{j \le m-1}
{\vert \nu \vert+j=m}} 
a_{\nu,j}(t) D^{\nu}_x D^{j}_t u=0, \quad t \ne 0,
\label{Equation}
\end{equation}
with the initial condition 
\begin{equation}
D^{k}_t u(0,x)=f_k(x) \in C^{\infty}_0(\mathbb{R}^n), \quad 
k=0,1,\cdots,m-1, \quad x \in \mathbb{R}^n, 
\label{Initial condition}
\end{equation}
where $D_t=-i\pa_t$ and $D^{\nu}_x=(-i\pa_{x_1})^{\nu_1}
\cdots (-i\pa_{x_n})^{\nu_n}$, 
$i=\sqrt{-1}$, for $\nu=(\nu_1,\ldots,\nu_n)$. 
Denoting by $\mathrm{Lip}_{\mathrm{loc}}(\R)$ 
the space of all functions which are locally Lipschitz on 
$\mathbb{R}$, we assume that each 
$a_{\nu,j}(t)$ belongs to $\mathrm{Lip}_{\mathrm{loc}}(\R)$ 
and satisfies 
\begin{equation}
a_{\nu,j}^\prime(t) \in L^1(\mathbb{R}) \quad 
\text{with $|\nu|+j=m$}.
\label{L1} 
\end{equation}
Moreover, following the standard definition of equations
of the regularly hyperbolic type
(e.g. Mizohata \cite{Mizohata}), we will assume
that the symbol of the differential operator 
$L(t,D_t,D_x)$ has real and distinct roots 
$\varphi_1(t;\xi),\ldots,\varphi_m(t;\xi)$ 
for $\xi \ne0$, and that
\begin{equation}
L(t,\tau,\xi)=(\tau-\varphi_1(t;\xi)) \cdots 
(\tau-\varphi_m(t;\xi)), 
\label{strict hyperbolicity1}
\end{equation} 
\begin{equation}
\inf_{\underset{j \ne k}
{\vert \xi \vert=1,t \in \mathbb{R}}} 
\vert \varphi_j(t;\xi)-\varphi_k(t;\xi) \vert>0.
\label{strict hyperbolicity2}
\end{equation}
Let us point out the main difficulties when trying to establish
dispersive estimates 
(i.e. the time decay estimates for the $L^p$--$L^q$ norms)
for equation \eqref{Equation}. Contrary
to the energy methods, for dispersive estimates
we need to have a good idea about the
propagators for the Cauchy problem 
\eqref{Equation}--\eqref{Initial condition}. Thus, we need
to make advances in the following two problems:
\begin{itemize}
\item to derive representation formulae for propagators for the
Cauchy problem \eqref{Equation}--\eqref{Initial condition} with
time-dependent coefficients. Ideally these propagators would
be in the form of oscillatory integrals;
\item to analyse the obtained representation formulae for
propagators taking into account the geometric properties of
characteristics which we know should be responsible for the
time decay rates of $L^p$--$L^q$ norms of propagators.
\end{itemize}
Thus, the aim of the paper is twofold. 
First, we will present representation formulae for propagators
for such equations. For this purpose we
will develop the asymptotic integration method
which is a parameter dependent version of the asymptotic
integration of ordinary differential equations 
(see e.g. \cite{Hartman}). 
We will trace the dependence on the parameter (which is
the frequency in this case) 
which is essential for further investigation.
This method, however, will present somewhat surprising
results. For example, the amplitudes of propagators
expressed in this form will have symbolic behaviour of
type $(0,0)$ rather than the usual $(1,0)$. Nevertheless,
this will be enough to carry out the second aim of this
part which is the further investigation
of the time decay properties of the propagators. 
The price that we will have to pay is that we may have
to assume additional regularity of the Cauchy data for high
frequencies. However, this is not so bad
because estimates for bounded times already will
require similar regularity assumptions. 

We will analyse the obtained representations to derive
time asymptotics of $L^p$--$L^q$ norms of the necessary
oscillatory integrals. There are several important differences
with the case of the wave equation, where level sets
of characteristics are nothing else but spheres, so one can simply
apply the stationary phase method to the obtained oscillatory
integrals. Now the critical points may be degenerate so the
stationary phase method (especially in the parameter depending
setting that we have here) does not work. 
In fact, the non-degeneracy of critical points is a rather
strong assumption for higher order equations, where degeneracy
of higher order may easily happen 
(examples of this are e.g. in \cite{sugi96}). That is why we
will allow them to be degenerate of a finite order, and the
time decay rates will depend on this order.
On the other hand,
van der Corput type estimates that are normally used in place
of the stationary phase in such 
problems are essentially one-dimensional, and
they do not take into account the geometric properties of 
characteristics (phases and characteristics 
do come from a hyperbolic equation after
all). So, we need to apply a parameter dependent version 
(now time is the parameter) of
van der Corput's lemma uniformly in $n-1$ directions of
non-vanishing higher order curvatures. Moreover, this has
to be done uniformly with respect to the time dependence
of the propagators. In fact,
we will relate the time-decay rates to the Sugimoto's
indices of levels sets of characteristic roots of the
limiting equation, thus establishing a more or less
complete picture of perturbation properties of dispersive
estimates for strictly hyperbolic equations with 
homogeneous symbols.
For example, in the
case of convex level sets one introduces the convex index 
$\gamma$ which is the largest order of tangency of tangent
lines to the level sets of characteristics of the limiting
equation and
it turns out that the $L^p$--$L^q$ norm of the 
corresponding propagator
decays as $t^{-\frac{n-1}{\gamma}}.$ In the case of the
second order equations one has $\gamma=2$ and
so one recovers that standard rate of decay of the wave
equation (see \cite{Bren75,Bren77,Litt73,stri70func}), 
and many other known results for the time independent
wave type second order equations. 
We also note that such index $\gamma$ does not play any
role for $L^p$--$L^p$ estimates, where singularities of
the projection from the canonical relations to the base
space start playing a role (see e.g. survey paper 
\cite{Ruzh-survey}).
The inclusion of mixed
derivatives in the symbols may influence the value of
$\gamma$. Moreover, mixed terms may make the analysis
more complicated. Already for the second order equations
this was demonstrated by Hirosawa and Reissig in \cite{HR},
for the problem of the influence of oscillations in coefficients.

In addition, methods introduced in this paper may be applied
to the study of strictly hyperbolic systems.
For example, let $A(t,D_x)$ be the first order $m\times m$ 
pseudo-differential system,
with entries $a_{ij}(t,\xi)$ being homogeneous with respect to
$\xi$ of order one and such that
$\pa_t a_{ij}(\cdot,\xi) \in L^1(\R)$ for all $\xi\in\Rn$.
We consider the evolution equation 
\begin{equation}\label{EQ:system1}
\pa_t U=iA(t,D_x)U,\quad U(0,x)=f(x),\; x\in\Rn.
\end{equation}
Let us assume that system \eqref{EQ:system1} is uniformly
strictly hyperbolic (see Mizohata \cite{Mizohata}), 
i.e. that its characteristics
$\va_k(t;\xi)$, $k=1,\cdots,m$, are real, and satisfy condition
\eqref{strict hyperbolicity2}. Then 
the strict hyperbolicity implies that 
we can diagonalise it similar
to Lemma \ref{lem:lem2.1}.
Thus, system \eqref{EQ:system1} splits into 
$m$ scalar first order equations of the form
$$
\partial_t v_k = i\va_k(t;D_x) v_k,\quad k=1,\ldots,m,
$$
for function $v_k$ related to the original vector function $U$. 
The condition on the integrability of time-derivatives of
$A$ implies that there is a limiting system
$A^\pm(t,\xi)=\lim_{t\to\pm\infty} A(t,\xi)$ with
characteristics $\va_k^\pm(\xi)=\lim_{t\to\pm\infty} \va_k
(t;\xi)$, which exist since we assume that
$\pa_t \va_k(\cdot;\xi)\in L^1(\R)$.
Therefore, solutions $v_k(t,x)$ can be analysed using estimates
for oscillatory integrals that we establish in \S 4. 
Details of this analysis are different from those for the
scalar equation \eqref{Equation}, especially in the
way of keeping track of the representation form for the
time derivatives of the solution, so we omit the analysis
of systems from this paper and it will appear elsewhere,
together with its specific applications and with refinements
of the analysis of high frequencies. We will also not
discuss the case of oscillations in this paper, but we
refer to, for example, the survey \cite{Reissig}, for the
overview of the case of the
wave equations. The case of oscillations
in higher order equations will appear elsewhere.
 
Thus, in \S 2 we will discuss the 
asymptotic integrations of the ordinary differential 
equations corresponding to 
our problem. Using these implicit representations, we will 
succeed to obtain 
the asymptotic integrations of 
\eqref{Equation}--\eqref{Initial condition}. 
The precise statement will be given in \S 3. 

Let us now give an informal overview of this method. 
Writing equation \eqref{Equation} as a system for
$$U=\left(|D|^{m-1}u,|D|^{m-2}D_t u,
\ldots,D^{m-1}_t u \right)^T$$
and taking 
the Fourier transform with respect to $x$, 
we can reduce it to the first order 
Cauchy problem 
\begin{equation}
D_t U=A(t;\xi) U,\quad U(0)=U_0.
\label{S1}
\end{equation} 
If we denote 
$$\vartheta_j(t;\xi)=\int^{t}_{0} \varphi_j(s;\xi) \, ds,\;
j=1,\ldots,m,$$ 
a natural candidate for the fundamental
matrix for \eqref{S1} is 
\[
\Phi(t;\xi)=\mathrm{diag}\left(e^{i\vartheta_1(t;\xi)}, 
\cdots,e^{i\vartheta_m(t;\xi)}\right).
\]
So, we look for the solution of \eqref{S1} in the form
$$U=\Phi(t;\xi) V, \;\; \textrm{with}\;\; 
V(t;\xi)=\pmb{\alpha}(\xi)+\pmb{\varepsilon}(t;\xi),$$
where we want $\pmb{\varepsilon}(t;\xi)$ to decay as $t\to\pm\infty$.
It can be checked that there is a matrix
$A_0(t;\xi)$ such that $D_t \Phi=A_0 \Phi$ and we get
$$
D_t U=D_t(\Phi V)=(D_t \Phi) V+\Phi D_t V=A_0 U+ \Phi D_t \pmb{\varepsilon}.
$$
Thus, $U$ becomes the solution of \eqref{S1} if we choose
$V$ and $\pmb{\varepsilon}$ such that
\begin{equation}
D_t V\equiv D_t \pmb{\varepsilon}=\Phi^{-1} (A-A_0) \Phi V.
\label{S2}
\end{equation} 
In \S 2 we will show that, in fact, there exists
a global-in-time solution $V$ of equation \eqref{S2}
of the required form $V=\pmb{\alpha}+\pmb{\varepsilon}$. Moreover,
$\pmb{\varepsilon}$ satisfies the property that
$\pmb{\varepsilon}(t;\xi)\to 0$ as $t\to\pm\infty$ for all $\xi\not=0$.
In addition, we will show the decay orders of both 
$\pmb{\alpha}(\cdot)$ and 
$\pmb{\varepsilon}(t;\cdot)$ 
and their derivatives. 
This will lead to an oscillatory integral representation
of solution $u(t,x)$ of \eqref{Equation} of the form
\begin{equation}
u(t,x)=\sum^{m-1}_{k=0}\sum^{m}_{j=1} \mathcal{F}^{-1} 
\left[ 
e^{i\vartheta_j(t;\xi)}
\left(\alpha^{j}_{k,\pm}(\xi)
+\varepsilon^{j}_{k,\pm}(t;\xi)\right)
 \widehat{f}_k (\xi)\right](x), 
\quad t \gtrless 0,
\label{S3}
\end{equation}
with amplitudes $\alpha^{j}_{k,\pm}(\xi), 
\varepsilon^{j}_{k,\pm}(t;\xi)$ of the form 
of $\pmb{\alpha}$ and 
$\pmb{\varepsilon}$ above. In fact, 
Theorem \ref{thm:asymptotic} will also
yield a similar representation for the derivatives of $u(t,x)$ 
with respect to time. The main difference with equations
with time independent
coefficients here is that the amplitudes $\alpha^{j}_{k,\pm}(\xi)$
and $\varepsilon^{j}_{k,\pm}(t;\xi)$ will have the symbolic
behavior of the type $(0,0)$ rather than the type $(1,0)$ usual
for equations with constant coefficients. 
Indeed, such choice of phases globally as
$\vartheta_{j}(t;\xi)$
introduces low order errors in the equation if we formally
substitute \eqref{S3} into \eqref{Equation} and as we know
the lower order terms may change the time decay properties in
an essential way (this is especially apparent for Schr\"odinger 
equations, but is also true in the hyperbolic case). Thus, the error
should be somehow accounted for and the behaviour of
amplitudes takes care of this. In any case, since we know that
the needed regularity of data comes from other parts of the
time-frequency phase space, we are still able to get the same time
decay rate under an additional regularity assumption in the
high frequency zone. So this difference does not change 
the final result in a big way.

Thus, in the second part of the paper we will use 
representation 
\eqref{S3} to derive the time decay of the $L^p$--$L^q$ 
norms of $u$, which in turn leads to Strichartz estimates 
and to well-posedness results for the corresponding 
semilinear equations in a rather (by now) standard way, so we
will derive the dispersive estimates and will
omit the details of the further standard analysis. 
In addition, in \S 4 we will present estimates for more general 
oscillatory integrals. Such estimates may be used not only
in the application to the problem we are considering in this
paper but in a wider range of applications. The estimates will
rely on estimates for parameter dependent oscillatory
integrals developed in \cite{R-vdC}. We state such result here
in Theorem \ref{THM:oscintthm}.
However, the meaning of the
parameter is different in our setting. 
Thus, in our problem here time
acts as a parameter while in problems for hyperbolic equations
with constant coefficients but with lower order terms considered
in \cite{RS2} the phase functions were not homogeneous and their
non-homogeneous contributions were considered to be a parameter
from the point of view of the perturbation theory. In principle,
it should be possible to combine problems with 
time-dependent coefficients  with those with lower order terms but
this will be a subject of another paper -- here we have an
advantage of making more use of the homogeneity of the
symbols and hence also of phases, considerably simplifying
some arguments. 
The obtained results can be applied to the global in time
well-posedness problems of Kirchhoff equations of high
orders and of Kirchhoff systems. Such applications
will be addressed elsewhere.


Let $\varphi_k^\pm(\xi)=\lim_{t\to\pm\infty}\varphi_k(t;\xi)$ 
be the limits of characteristic roots as will be shown 
to exist in 
Proposition \ref{prop:root}. Let us introduce the convex
and non-convex Sugimoto indices for the level sets of these
functions. In the time independent setting these indices have
been introduced by Sugimoto in \cite{sugi94,sugi96}.
These indices will determine the decay rate of propagators
for large frequencies. 

Let $\varphi\in C^\infty(\Rn\backslash 0)$ be a homogeneous of
order one function and let 
$\Sigma_\varphi=\{\xi\in\Rn:\; \varphi(\xi)=1\}$
be its level set. Suppose first that $\Sigma_\varphi$ is convex.
We define the {\em convex Sugimoto index}
$\gamma(\Sigma_\va)$ of $\Sigma_\va$ by 
\begin{equation}\label{eq:convex-sugimoto}
\gamma(\Sigma_\varphi):=\sup_{\sigma\in\Sigma_\varphi}\sup_P\gamma
(\Sigma_\varphi;\sigma,P)\,,
\end{equation}
where $P$ is a plane containing the normal to~$\Sigma_\varphi$ at~
$\sigma$ and
$\gamma(\Sigma_\varphi;\sigma,P)$ denotes the order of the contact 
between the line $T_\sigma\cap P$ (where $T_\sigma$ is the tangent 
plane at~$\sigma$), and the curve $\Sigma_\varphi\cap P$. 

In the case when the level set $\Sigma_\varphi$ is not convex, we get a 
weaker result based on the van der Corput lemma. In this case we
use the {\em non-convex Sugimoto index} 
$\gamma_0(\Sigma_\va)$ of $\Sigma_\varphi$ which
we define as
\begin{equation}\label{eq:nonconvex-Sugimoto}
\gamma_0(\Sigma_\varphi):=\sup_{\sigma\in\Sigma_\varphi}\inf_P\gamma
(\Sigma_\varphi;\sigma,P)\,,
\end{equation}
where $P$ and $\sigma$ are the same as in the convex case.

We note that for the second order equations we have
$\gamma=\gamma_0=2$ and the following theorem covers the case of the
wave equation as a special case, also improving the 
corresponding result in \cite{Matsuyama1}. 
We use the notation $L^{p}_s(\Rn)$ for the
standard Sobolev space with $s$ derivatives over 
$L^p(\mathbb{R}^n)$, and by $\dot{L}^p_s(\Rn)$ we denote its
homogeneous version.
The result on the dispersive
estimates that we will prove among other things, is as follows:

\begin{thm}\label{thm:Main theorem}
Assume \eqref{L1}--\eqref{strict hyperbolicity2}. 
Then the solution $u(t,x)$ of \eqref{Equation} 
satisfies the following 
estimates{\rm :}

{\rm (i)} Suppose that the set
$$\Si_{\varphi^{\pm}_k}=
\{\xi\in\Rn: \va^{\pm}_k(\xi)=1\}$$ 
is convex for all $k=1,\ldots,m$, and set
$\gamma=\underset{k=1,\ldots,m}{\max} 
\gamma(\Sigma_{\varphi^{\pm}_k})$. 
In addition, suppose that 
$(1+\vert t \vert)^r a^\prime_{\nu,j} 
\in L^1(\mathbb{R})$ for 
$1\le r \le [(n-1)/\gamma]+1$, 
and for all $\nu,j$ with $|\nu|+j=m$. 
Let $1<p\le 2 \le q<+\infty$ and $\frac{1}{p}+\frac{1}{q}=1$. 
Then for all $t\in \mathbb{R}$ we have the estimate
\begin{multline}\label{EQ:mainest1}
\Vert D^{l}_t D^{\alpha}_x u(t,\cdot) \Vert_{L^q(\Rn)}
\leq \\
C (1+|t|)^{-\frac{n-1}{\gamma}\p{\frac{1}{p}- \frac{1}{q}}}
\sum^{m-1}_{k=0} \left(\|f_k \|_
{\dot{L}^p_{N_p+l+\vert \alpha\vert-k}(\Rn)}+
\|f_k \|_{\dot{L}^{p}_{l+|\alpha|-k}(\Rn)}\right),
\end{multline}
where $N_p=\left(n-\frac{n-1}{\gamma}+
\left[ \frac{n-1}{\gamma} \right]+1\right)
\left(\frac{1}{p}-\frac{1}{q}\right)$, $l=0,\ldots,m-1$, and 
$\alpha$ is any multi-index.

{\rm (ii)} Suppose that $\Si_{\varphi^{\pm}_k}$ 
is non-convex 
for some 
$k=1,\ldots,m$, and let us set $\gamma_0=
\underset{k=1,\ldots,m}{\max} 
\gamma_0(\Sigma_{\varphi^{\pm}_k})$. 
In addition, suppose that 
$(1+\vert t \vert) a^\prime_{\nu,j} \in L^1(\mathbb{R})$ 
for all $\nu,j$ with $|\nu|+j=m$. Let $1<p\le 2 \le q<+\infty$ 
and $\frac{1}{p}+\frac{1}{q}=1$. 
Then for all $t\in \mathbb{R}$ we have
the estimate 
\begin{multline*}
\| D^{l}_t D^{\alpha}_x u(t,\cdot) \Vert_{L^q(\Rn)} 
\leq \\
C (1+|t|)^{-\frac{1}{\gamma_0}\p{\frac{1}{p}- \frac{1}{q}}}
\sum^{m-1}_{k=0} \left(\|f_k \|_
{\dot{L}^p_{N_p+l+\vert \alpha\vert-k}(\Rn)}+
\|f_k \|_{\dot{L}^{p}_{l+|\alpha|-k}(\Rn)}\right),
\end{multline*}
where $N_p=\left(n-\frac{1}{\gamma_0}+1\right)
\left(\frac{1}{p}-\frac{1}{q}\right)$, $l=0,\ldots,m-1$, and 
$\alpha$ is any multi-index.
\end{thm}

\begin{rem}
The way we formulate the estimates in Theorem
\ref{thm:Main theorem} is to unify different estimates for
different parts of the solution. 
This may explain the appearance of two norms in the right
hand side of \eqref{EQ:mainest1}, for example, to account
for both small and large frequencies.
The much more precise
estimates are possible and they are stated in Theorem
\ref{thm:Main theorem-2}. 
\end{rem}

Let us now make only a few short
remarks to compare our results with what is 
known for $m=2$. 
In the constant coefficient case and $p=2$
the estimates coincide with those for constant coefficient
equations considered in \cite{sugi94}. 
Also, in the case of constant coefficients, the Sobolev
index $N_p$ in the 
estimate \eqref{EQ:mainest1} can be improved since there
is no addition of the integer part in its definition in this
case. For the detailed overview of constant 
coefficients case we refer to \cite{RS2}, 
and results in this direction for non-constant
coefficients were announced in
\cite{MR09} and the detailed proofs will appear elsewhere.

Further, in the time-dependent case, 
it can be already noted that the statement of Theorem
\ref{thm:Main theorem} goes beyond results
available in certain energy classes. For example, in the
often considered case of the time--dependent wave equation
(so that $m=2$, 
e.g. \cite{HR}, \cite{Matsuyama2}, \cite{Reissig1},
\cite{Reissig}, etc.) one
obtains the estimate for
$\Vert D_t u(t,\cdot) \Vert_{L^q(\Rn)}$ and 
$\Vert \nabla u(t,\cdot) \Vert_{L^q(\Rn)}$ only,
and not for solution $u$ itself. Moreover,
the use of homogeneous Sobolev spaces in 
\eqref{EQ:mainest1} allows to gain more
information in the low frequency region. 
At the same time, also already 
for the case $m=2$, we make the assumption on only 
one derivative of the coefficients $a_{\nu,j}$, which is
another improvement compared with the known literature.
This improvement will be crucial in dealing with applications
to Kirchhoff equations.

We will denote $\jp{x}=\sqrt{1+|x|^2}$. Constants may change
from formula to formula, although they are usually denoted by
the same letter.

The authors thank Jens Wirth for remarks leading to an
improvement of the manuscript.

\section{Asymptotic integration of ODE} 
\setcounter{equation}{0} 
In this section we will construct the asymptotic integration 
of the ordinary differential equation.
By applying the Fourier transform on $\mathbb{R}^n_x$ 
to \eqref{Equation}, 
we get 
\begin{equation}
D^{m}_t v+\sum^{m}_{j=1} h_j(t;\xi)D^{m-j}_{t} v=0, 
\label{wintner1}
\end{equation} 
where 
\[
h_j(t;\xi)=\sum_{\vert \nu \vert=j} a_{\nu,m-j}(t) \xi^\nu, 
\quad \xi \in \mathbb{R}^n
\]
(note that there is a slight change of the meaning of $j$ here
compared to \eqref{Equation}). 
This is the ordinary differential equation of homogeneous 
$m^{\rm th}$ order 
with the parameter $\xi=(\xi_1,\ldots,\xi_n)$. As usual, the strict 
hyperbolicity 
\eqref{strict hyperbolicity1}--\eqref{strict hyperbolicity2} 
means that the characteristic roots of \eqref{wintner1} 
are real and distinct. We denote them by 
$\varphi_1(t;\xi),\ldots,\varphi_m(t;\xi)$. Notice that each 
$\varphi_\ell(t;\xi)$ has a homogeneous degree of order one
with respect to 
$\xi$. In this section we will develop an asymptotic integration of 
the equation \eqref{wintner1}. Let us start by writing 
\eqref{wintner1} as the first 
order system. In \eqref{wintner1} we put for brevity
\[
H_j(t,\xi)=h_j(t;\xi/\vert \xi \vert), 
\]
and denote 
\[
v_j=\vert \xi \vert^{m-j-1} D^{j}_t v, \quad j=0,1,\cdots,m-1. 
\]
It is easy to see that 
\[
D_t v_j=\vert \xi \vert v_{j+1}, \quad j=0,1,\cdots,m-2,
\]
holds. Then \eqref{wintner1} can be written as 
\[
D_t v_{m-1}+\sum^{m-1}_{j=0} H_{m-j}(t,\xi) 
\vert \xi\vert v_j=0. 
\]
Hence, if we put 
\[
\mathcal{H}(t,\xi)=\left(
\begin{array}{cccc}
0 & 1 & \ldots & 0 \\ 
0 & 0 & \ddots & 0 \\
\vdots & \ddots & \ddots & 1 \\
-H_m(t,\xi) & -H_{m-1}(t,\xi) & \ldots & -H_1(t,\xi) 
\end{array}\right), 
\]
then \eqref{wintner1} can be written by 
\begin{equation}
D_t \pmb{v}=\mathcal{H}(t,\xi) \vert \xi \vert \pmb{v},
\label{system}
\end{equation}
where $\pmb{v}={}^T (v_0,v_1,\ldots,v_{m-1})$.

We will use the following lemma. 
\begin{lem}[\cite{Mizohata}~Proposition~6.4] \label{lem:lem2.1}
Assume \eqref{L1}--\eqref{strict hyperbolicity2}. Then 
there exists a matrix 
$\mathcal{N}=\mathcal{N}(t;\xi)$ of homogeneous 
degree $0$ satisfying the following properties{\rm :} 

\noindent 
{\rm (i)} $\mathcal{N}\mathcal{H}(t,\xi)=
\mathcal{D}\mathcal{N}$, where 
\[\mathcal{D}=\mathcal{D}(t;\xi)=
\mathrm{diag} \left\{ \varphi_1(t;\xi/|\xi|),\ldots,
\varphi_m(t;\xi/|\xi|)\right\},
\]

\noindent 
{\rm (ii)} $\displaystyle
{\inf_{\xi \in \mathbb{R}^n\backslash 0, t \in \mathbb{R}}}
\vert {\rm det} \, \mathcal{N}(t;\xi) \vert >0$, 

\noindent 
{\rm (iii)} $\mathcal{N}(t;\xi)$ is $C^\infty$ in 
$\xi \ne 0$, $C^1$ in $t \in \mathbb{R}$ and 
$\partial_t \mathcal{N}(t;\xi)$ belongs 
to $L^1(\mathbb{R})$ for each $\xi\not=0$. 
\end{lem}

We will first derive the energy estimates. 
\begin{lem} \label{lem:lem2.2}
Assume \eqref{L1}--\eqref{strict hyperbolicity2}. Let 
$\pmb{v}=\pmb{v}(t;\xi)$ be a general solution of \eqref{system}. 
Then, for all $t \in \mathbb{R}$, we have 
\begin{equation}\label{energy}
\vert \pmb{v}(t;\xi) \vert^2 \le C \vert \pmb{v}(0;\xi) \vert^2 
\mathrm{e}^{\int^{+\infty}_{-\infty}
2\Vert \partial_t \mathcal{N}(s;\xi) \Vert \, ds}. 
\end{equation}
\end{lem} 
\begin{proof} Multiplying \eqref{system} by $\mathcal{N}=\mathcal{N}(t;\xi)$ 
from Lemma~\ref{lem:lem2.1}, we get 
\[
D_t (\mathcal{N}\pmb{v})-\mathcal{N}\mathcal{H}\vert \xi \vert \pmb{v}
-(D_t \mathcal{N}) \pmb{v}=0.
\]
Putting $\mathcal{N} \pmb{v}=\pmb{w}$, we have 
\begin{equation}\label{EQ:Origine}
D_t \pmb{w}-\mathcal{D}\vert \xi \vert \pmb{w}-(D_t \mathcal{N})\pmb{v}=0, 
\end{equation}
since $\mathcal{N}\mathcal{H}=\mathcal{D}\mathcal{N}$ by 
Lemma \ref{lem:lem2.1}. This implies that 
\[
\partial_t \vert \pmb{w} \vert^2 
= 2 \Re \left( \partial_t\pmb{w} \cdot \overline{\pmb{w}} \right) 
= 2 \Re \left( i\mathcal{D}\vert \xi \vert \pmb{w} \cdot 
\overline{\pmb{w}}\right) 
+2\Re \left(i (D_t \mathcal{N}) \pmb{v} \cdot \overline{\pmb{w}} \right).
\]
Taking account that $\overline{i\mathcal{D}\vert \xi \vert}=-i\mathcal{D}
\vert \xi \vert$ and $\mathcal{D}$ is real and diagonal, we have 
$$\Re \left( i\mathcal{D}\vert \xi \vert \pmb{w} \cdot 
\overline{\pmb{w}}\right)=0,
$$ 
hence, 
\begin{equation}
\partial_t \vert \pmb{w} \vert^2 \le 2\Vert \partial_t \mathcal{N} \Vert 
\vert \pmb{v} \vert \vert \pmb{w} \vert. 
\label{differential inequality}
\end{equation}
Here, $\vert \pmb{v} \vert$ and $\vert \pmb{w} \vert$ are equivalent to each 
other. Indeed, there exists $C_1,C_2>0$ such that $C_1 \vert \pmb{v} \vert \le 
\vert \pmb{w} \vert \le C_2 \vert \pmb{v} \vert$ on account of 
Lemma \ref{lem:lem2.1}. Thus, integrating \eqref{differential inequality}, we 
arrive at 
\[
\vert \pmb{v}(t;\xi) \vert^2 \le C \left( \vert \pmb{v}(0;\xi) \vert^2 
+\int^{\vert t \vert}_{0} \Vert \partial_t \mathcal{N}(s;\xi) \Vert 
\vert \pmb{v}(s;\xi) \vert^2 \, ds \right). 
\]
Since $\partial_t \mathcal{N} \in L^1(\mathbb{R})$ by Lemma \ref{lem:lem2.1} 
(iii), we conclude from Gronwall's lemma that \eqref{energy} is true. 
The proof of Lemma \ref{lem:lem2.2} is complete. 
\end{proof} 

As a consequence of \eqref{EQ:Origine} in the proof of 
Lemma \ref{lem:lem2.2}, we have derived 
\begin{equation}\label{EQ:Orignal equation}
D_t \pmb{w}-\mathcal{D}|\xi| \pmb{w}
-(D_t \mathcal{N})\mathcal{N}^{-1}\pmb{w}=0,
\end{equation}
where we put $\pmb{w}=\mathcal{N}\pmb{v}$. 
We can expect that the solution of \eqref{EQ:Orignal equation} 
is asymptotic to the solution of 
\begin{equation}\label{EQ:Compared equation}
D_t \pmb{y}=\mathcal{D}|\xi|\pmb{y}.
\end{equation}
Let $\Phi(t;\xi)$ be the fundamental matrix of \eqref{EQ:Compared equation}, 
i.e., 
\[
\Phi(t;\xi)=
\mathrm{diag} \left\{ 
e^{i\vartheta_1(t;\xi)}, \cdots, e^{i\vartheta_m(t;\xi)}
\right\}, 
\]
where we put 
\[
\vartheta_j(t;\xi)=\int_0^t \varphi_j(s;\xi)\, ds, \qquad 
j=1,\ldots,m.
\]

Let us first analyse certain basic properties of characteristic roots 
$\varphi_k(t;\xi)$ of \eqref{strict hyperbolicity1}.  

\begin{prop} \label{prop:root}
Let the operator $L(t,D_t,D_x)$ satisfy the properties
\eqref{strict hyperbolicity1}--\eqref{strict hyperbolicity2}. Then each 
$\partial_t \varphi_k(t;\xi)$, $k=1,\ldots,m$, is 
homogeneous of order one in $\xi$, and there exist a constant $C>0$ such that 
\begin{equation}\label{EQ:phi-bounded0}
| \partial_t\varphi_k(t;\xi)|\leq C|\xi| 
\quad \textrm{for all}\;\; \xi\in\Rn, 
\;\; t \in \mathbb{R}, \;\; k=1,\ldots,m.
\end{equation}
Moreover, if $a^\prime_{\nu,j}(\cdot)\in L^1(\R)$ for all $\nu,j$, 
then we have also 
$\partial_t\varphi_k(\cdot;\xi)\in L^1(\R)$ for all 
$\xi\in\mathbb{R}^n$. Furthermore, there exist functions 
$\varphi_k^\pm\in C^\infty(\Rn\backslash 0)$, homogeneous of order one, 
such that 
\begin{equation}\label{EQ:convphi}
\partial_\xi^\alpha\varphi_k(t;\xi)\to 
\partial_\xi^\alpha\varphi_k^\pm(\xi) \textrm{ as } 
t\to\pm\infty,
\end{equation} 
for all $\xi\in\Rn$, all $\alpha$, and $k=1,\ldots,m$. 
Finally, we have the following 
formula for the derivatives of characteristic roots{\rm :}
\begin{equation}
\partial_t\varphi_k(t;\xi)= 
-\sum_{|\nu|+j=m} a^\prime_{\nu,j}(t) \varphi_k(t;\xi)^j\xi^\nu
\prod_{r\not=k} 
\p{\varphi_k(t;\xi)-\varphi_r(t;\xi)}^{-1}.
\label{product0}
\end{equation}
\end{prop}
\begin{proof}
Let us show first that $\varphi_k(t;\xi)$ is 
bounded with respect to 
$t\in \mathbb{R}$, i.e., 
\begin{equation}
\vert \varphi_k(t;\xi) \vert \le C \vert \xi \vert, \quad 
\text{for all $\xi \in \mathbb{R}^n$, $t \in \mathbb{R}$, $k=1,\ldots,m$.}
\label{phibdd}
\end{equation}
We will use the fact that $\varphi_k(t;\xi)$ 
are roots of the polynomial 
$L$ of the form 
$$L(t,\tau,\xi)=\tau^m+c_1(t,\xi)\tau^{m-1}
+\cdots+c_m(t,\xi)$$ 
with $|c_j(t,\xi)|\leq M|\xi|^j$, for some $M\geq 1$. 
Suppose that one of its roots $\tau$ satisfies 
$|\tau(t,\xi)|\geq 2M|\xi|$. 
Then
\begin{align*}|L(t,\tau,\xi)| & 
 \geq |\tau|^{m} \p{1-\frac{|c_1(t,\xi|)}{|\tau|}-
\cdots-\frac{|c_m(t,\xi)|}{|\tau|^m} }      \\
& \geq 
2M|\xi|^m \p{1-\frac{1}{2}-\frac{1}{4M}-\cdots-\frac{1}{2^m M^{m-1}}}>0,
\end{align*}
hence $|\tau(t,\xi)|\leq 2M|\xi|$ for all $\xi\in\Rn$. 
Thus we establish \eqref{phibdd}. 

Differentiating \eqref{strict hyperbolicity1} 
with respect to $t$, we get
$$ 
\frac{\partial L(t,\tau,\xi)}{\partial t}=
\sum_{|\nu|+j=m} a^\prime_{\nu,j}(t) \tau^j\xi^\nu=
-\sum_{k=1}^m \partial_t\varphi_k(t;\xi)\prod_{r\not=k} 
\p{\tau-\varphi_r(t;\xi)}.
$$
Setting $\tau=\varphi_k(t;\xi)$, we obtain
\begin{equation}
\partial_t\varphi_k(t;\xi)\prod_{r\not=k} 
\p{\varphi_k(t;\xi)-\varphi_r(t;\xi)}= -
\sum_{|\nu|+j=m} a^\prime_{\nu,j}(t) \varphi_k(t;\xi)^j\xi^\nu,
\label{product}
\end{equation}
implying \eqref{product0}.
Now, using \eqref{strict hyperbolicity2}, \eqref{phibdd}, 
and the assumption 
that $a^\prime_{\nu,j}(\cdot)\in L^1(\R)$ for all $\nu, j$, 
we conclude that \eqref{EQ:phi-bounded0} holds and 
$\partial_t\varphi_k(\cdot;\xi)\in L^1(\R)$ 
for all $\xi\in\Rn$ and 
$k=1,\ldots,m$. 
The homogeneity of order one of $\partial_t \varphi_k(t;\xi)$ 
is an immediate consequence of \eqref{product} 
and its derivatives. 

Finally, setting 
$$\varphi_\ell^\pm(\xi)=\varphi_k(0;\xi)
+\int_0^{\pm\infty}\partial_t\varphi_k(t;\xi) \, dt,$$ 
we get \eqref{EQ:convphi} with $\alpha=0$. Differentiating
this equality with respect to $\xi$, we get
\eqref{EQ:convphi} for all $\alpha$. The proof is complete.
\end{proof}

We note that under the assumptions of 
Proposition \ref{prop:root} the 
coefficients $a_{\nu,j}(t)$ of the operator $L(t,D_t,D_x)$ 
in \eqref{Equation} 
have limits $a_{\nu,j}^\pm$ as $t\to\pm\infty$, namely 
$$a_{\nu,j}^\pm=a_{\nu,j}(0)+\int_{0}^
{\pm\infty} a^\prime_{\nu,j}(t) \, dt.$$ 
Functions $\varphi^\pm_k(\xi)$ are characteristics of the 
limiting strictly hyperbolic operator
\begin{equation}
L^\pm(D_t,D_x)u 
\equiv D^{m}_t u+\sum_{\underset{j \le m-1}
{\vert \nu \vert+j=m}} 
a^\pm_{\nu,j} D^{\nu}_x D^{j}_t u,
\label{Equation-lim}
\end{equation}
and their geometric properties are responsible for the time decay
of solutions to the Cauchy problems for both operators
$L(t,D_t,D_x)$ and $L^\pm(D_t,D_x)$. 
This will be analysed in \S 4. 
We also note that since operator $L^\pm(D_t,D_x)$
has constant coefficients its solution can be represented
as a sum of oscillatory integrals in the standard way.
The dependence of coefficients of $L(t,D_t,D_x)$ on time brings
corrections to the phases and amplitudes of this representation.

Next we make the representation formulae of 
solutions for our equation. 
The following proposition is known as Levinson's lemma 
(see Coddington and Levinson \cite{Coddington}) in the 
theory of ordinary differential equations, the new feature
here is the additional dependence on $\xi$. 
For the convenience of the readers, we shall prove 
it along the method of 
Ascoli \cite{Ascoli} and Wintner \cite{Wintner} 
(cf. Hartman \cite{Hartman}).

\begin{prop} \label{prop:prop2.3} 
Assume \eqref{L1}--\eqref{strict hyperbolicity2}. Then, for 
every nontrivial solution $v(t;\xi)$ of \eqref{wintner1}, there exist 
vectors of $C^{\infty}$-amplitude functions 
$\pmb{\alpha}_\pm(\xi)$ and error functions 
$\pmb{\varepsilon}_\pm(t;\xi)$ such 
that 
\begin{equation}\label{converse1} 
\pmb{v}(t;\xi)=\mathcal{N}(t;\xi)^{-1}\Phi(t;\xi) 
\left( \pmb{\alpha}_\pm(\xi)
+\pmb{\varepsilon}_\pm(t;\xi) \right), \quad t \gtrless 0, 
\end{equation} 
where 
\begin{equation} 
\pmb{\varepsilon}_\pm(t;\xi) \to \pmb{0} 
\quad \text{for any fixed $\xi \ne 0$ as $t \to \pm\infty$.} 
\label{converse3} 
\end{equation}
Furthermore we have 
\begin{equation}
\label{add}
D_t\pmb{\varepsilon}_\pm(t;\xi)
=C(t;\xi)\left(\pmb{\alpha}_\pm(\xi)+ \pmb{\varepsilon}_\pm(t;\xi) 
\right),  
\end{equation}
where $C(t;\xi)$ belongs 
to $L^1(\mathbb{R})$ in $t$, and 
has the following form{\rm :} 
\begin{equation}
\label{matrix norm}
C(t;\xi)=\Phi(t;\xi)^{-1} (D_t \mathcal{N}(t;\xi))
\mathcal{N}(t;\xi)^{-1}\Phi(t;\xi).
\end{equation}
\end{prop}
\begin{proof} We can expect that every solution $\pmb{w}=\pmb{w}(t;\xi)$ of 
\eqref{EQ:Orignal equation} is asymptotic to some solution 
$\pmb{y}=\pmb{y}(t;\xi)$ of \eqref{EQ:Compared equation}. 
If we perform the Wronskian transform $\pmb{z}=\Phi(t;\xi)^{-1}\pmb{w}$, 
then the system \eqref{EQ:Orignal equation} reduces to a system 
$D_t \pmb{z}=C(t;\xi) \pmb{z}$, where $C(t;\xi)$ is given by 
\eqref{matrix norm}. 

We will now prove that \eqref{converse1}--\eqref{converse3} 
hold for every nontrivial 
solution $v=v(t;\xi)$. It follows from Lemma \ref{lem:lem2.2} that 
\[
|\pmb{z}(t;\xi)| \le 
\|\Phi(t;\xi)^{-1}\| |\pmb{w}(t;\xi)| \le 
c |\pmb{w}(t;\xi)| \le c_1 |\pmb{w}(0;\xi)|
\]
for all $t \in \mathbb{R}$ and some constant $c_1$. 
Using this bound and 
equation $D_t \pmb{z}=C(t;\xi)\pmb{z}$, we have 
\[
|D_t \pmb{z}(t;\xi)| \le \|C(t;\xi)\|\, |\pmb{z}(t;\xi)|
\le c_1 |\pmb{w}(0;\xi)| \| C(t;\xi) \|, 
\]
hence $D_t \pmb{z}(\cdot;\xi) \in L^1(\mathbb{R})$ on account 
of $C(\cdot;\xi)\in 
L^1(\mathbb{R})$. Thus $\{ \pmb{z}(t;\xi)\}_{t \in \mathbb{R}}$ 
is a convergent function, and there exists 
$$\displaystyle{\lim_{t \to \pm\infty}}
\pmb{z}(t;\xi)=:\pmb{\alpha}_{\pm}(\xi).
$$ 
If we set 
\begin{equation}\label{EQ:Z-EQ}
\pmb{\varepsilon}_\pm(t;\xi)=\pmb{z}(\xi,t)-\pmb{\alpha}_{\pm}(\xi),
\end{equation}
then $\pmb{z}(t;\xi)$ can be written 
as $\pmb{z}(t;\xi)=\pmb{\alpha}_\pm(\xi)+
\pmb{\varepsilon}_\pm(t;\xi)$ 
for $t \gtrless 0$, and further, $\pmb{\varepsilon}_\pm(t;\xi)$ decays 
as $t \to \pm\infty$ for any fixed $\xi \ne 0$, which proves \eqref{converse3}. Since $\pmb{w}(t;\xi)=\Phi(t;\xi) \pmb{z}(t;\xi)$, we get the formula 
\eqref{converse1}. Finally, differentiating \eqref{EQ:Z-EQ} with 
respect to $t$ and using the equation $D_t \pmb{z}=C(t;\xi) \pmb{z}$, we get 
\eqref{add}. The proof of Proposition 
\ref{prop:prop2.3} is now complete. 
\end{proof}  

Finally, we will need the estimates for higher 
order derivatives of $C(t;\xi)$ 
appearing in Proposition \ref{prop:prop2.3}.
\begin{lem} \label{lem:derivative of C} 
Assume \eqref{L1}--\eqref{strict hyperbolicity2}. 
Then the $\mu^{\rm th}$ 
derivatives for each entry $c_{jk}(t;\xi)$ of $C(t;\xi)$ satisfy
\begin{equation} 
\left| \partial^{\mu}_{\xi} c_{jk}(t;\xi) \right|
\le c |\xi|^{-|\mu|}
(1+|t|)^{|\mu|} \Psi(t),
\quad j,k=1,\ldots,m, 
\label{growth}
\end{equation}
for $|\mu|\ge 1$ and $|\xi|\ge 1$, where 
\begin{equation}\label{Psi}
\Psi(t)=\sum_{\underset{j\le m-1}{|\nu|+j=m}} 
\left| a^\prime_{\nu,j}(t)\right|. 
\end{equation}
For $0<|\xi|<1$ and $|\mu|\ge 1$, we have 
\begin{equation} 
\left|\partial^{\mu}_{\xi} c_{jk}(t;\xi) \right|
\le c |\xi|^{-|\mu|}
(1+|t|)^{|\mu|} \Psi(t),
\quad j,k=1,\ldots,m. 
\label{low growth}
\end{equation}
Moreover, assume that 
$(1+|t|)^{|\mu|}a^\prime_{\nu,j}(t) 
\in L^1(\mathbb{R})$ for all $\nu,j$, and for some 
$\mu$ with $|\mu| \ge 1$. Then we have 
$\partial^{\mu}_{\xi} c_{jk}(t;\xi)\in L^1(\mathbb{R})$. 
\end{lem}
\begin{proof}
Since $\varphi_j(t;\xi)$ is homogeneous of order one, we have 
\[
\vert \partial^{\mu}_\xi \varphi_j(t;\xi) \vert 
\le c \vert \xi \vert^{-\vert \mu \vert+1} 
\quad \text{for $\xi\in \mathbb{R}^n\setminus 0$, $j=1,\ldots,m$,} 
\]
and hence, 
\[
\vert \partial^{\mu}_\xi \vartheta_j(t;\xi)\vert \le 
\int^{|t|}_{0} \vert \partial^{\mu}_\xi \varphi_j(s;\xi)\vert \, ds 
\le c|t| |\xi|^{-\vert \mu \vert+1}, 
\quad \text{for $\xi\in \mathbb{R}^n\setminus 0$, $j=1,\ldots,m$.} 
\]
Thus we get, for every multi-index $\mu$, 
\begin{equation}
\left\vert \partial^{\mu}_\xi e^{i\vartheta_j(t;\xi)} \right\vert 
\le \begin{cases}
c(1+\vert t \vert)^{|\mu|}, & \text{$|\xi|\ge 1$,}\\ 
c(1+\vert t \vert)^{|\mu|}|\xi|^{-|\mu|+1}, & \text{$0<|\xi|<1$.}
\end{cases}
\label{theta}
\end{equation}
Now let us go back to \eqref{matrix norm}. It follows from 
\eqref{product0} that $\partial_t \mathcal{N}$ is represented by 
$a^\prime_{\nu,j}$. 
Hence, using \eqref{theta} and differentiating \eqref{matrix norm} with 
respect to $\xi$, we conclude that $\mu^{\rm th}$ derivative of 
$c_{jk}(t;\xi)$ with respect to $\xi$ is bounded by 
$|\xi|^{-|\mu|}
(1+|t|)^{|\mu|} \Psi(t)$ for 
$|\xi| \ge 1$, where $\Psi(t)$ is given in \eqref{Psi}. 
$\partial^{\mu}_{\xi} c_{jk}(t;\xi)\in L^1(\mathbb{R})$ follows from 
the assumption 
$(1+|t|)^{|\mu|}a^\prime_{\nu,j}(t)
\in L^1(\mathbb{R})$ for all $\nu,j$. The proof of 
Lemma \ref{lem:derivative of C} is complete. 
\end{proof}

\section{Representation of solution} 
In this section we will establish the representation formulae
for solutions of the Cauchy problem \eqref{wintner1} in the
form of the oscillatory integrals. Let $$V(t;\xi)=
{}^T (\pmb{v}_0(t;\xi),\ldots,\pmb{v}_{m-1}(t;\xi))$$ 
be the fundamental matrix 
of \eqref{EQ:Orignal equation}. This means that $V(0;\xi)=I$. 
Then it follows from Proposition \ref{prop:prop2.3} that 
each $\pmb{v}_j(t;\xi)$ can be 
represented by 
\begin{equation}
\pmb{v}_j(t;\xi)=\mathcal{N}(t,\xi)^{-1}\Phi(t;\xi)
\left(\pmb{\alpha}_{j,\pm}(\xi)+
\pmb{\varepsilon}_{j,\pm}(t;\xi)\right). 
\label{Base}
\end{equation}
Let $u(t,x)$ be the solution to \eqref{Equation} 
with the Cauchy data 
$D^{k}_t u(0,x)=f_k(x)$. Put 
$$\widehat{\pmb{u}}(t,\xi)=
\left(|\xi|^{m-1}\widehat{u}(t,\xi),\ldots,|\xi|^{m-1-l}D_t^l 
\widehat{u}(t,\xi),
\ldots,D_t^{m-1}\widehat{u}(t,\xi) \right)^T.$$ 
Then we can write 
$\widehat{\pmb{u}}(t,\xi)=V(t;\xi)\widehat{\pmb{u}}(0,\xi)$; 
thus we arrive at 
\begin{equation}
|\xi|^{m-1-l}D_t^l 
\widehat{u}(t,\xi)=\sum_{j=1}^m \sum_{k=0}^{m-1}  e^{i\vartheta_j(t;\xi)}
n^{l j}(t;\xi)
\left(\alpha_{j,\pm}^k(\xi)+\varepsilon_{j,\pm}^k(t;\xi)\right)
|\xi|^{m-1-k}\widehat{f}_k(\xi),
\end{equation}
where $n^{l j}(t;\xi)$ is the entry of $\mathcal{N}(t;\xi)^{-1}$: 
\[
\mathcal{N}(t;\xi)^{-1}
=\left(n^{l j}(t;\xi)\right)_{\underset{j=1,\ldots,m}{l=0,\ldots,m-1}}.
\]

Summarizing the above argument, 
we have the representation formulae of 
\eqref{Equation}--\eqref{Initial condition}. 
\begin{thm} \label{thm:asymptotic}
Assume \eqref{L1}--\eqref{strict hyperbolicity2}. 
Then there exists 
$\alpha^{j}_{k,\pm}(\xi)$ and 
$\varepsilon^{j}_{k,\pm}(t;\xi)$ such that 
the solution 
$u(t,x)$ of our problem 
$\eqref{Equation}$--$\eqref{Initial condition}$ is 
represented by 
\[
D_t^l u(t,x)=\sum_{j=1}^m \sum_{k=0}^{m-1} \mathcal{F}^{-1}\left[ 
e^{i\vartheta_j(t;\xi)}
n^{l j}(t;\xi)
\left(\alpha_{j,\pm}^k(\xi)+\varepsilon_{j,\pm}^k(t;\xi)\right)
|\xi|^{l-k}\widehat{f}_k(\xi) \right](x),
\quad t \gtrless 0,
\] 
for $l =0,\ldots,m-1$, where 
\[
\left| \alpha^k_{j,\pm}(\xi)\right|\le c, 
\quad \left| \varepsilon^k_{j,\pm}(t;\xi)\right| \le c
\int^{+\infty}_{\vert t \vert} \Psi(s) \, ds, 
\] 
and $\Psi(t)$ is defined by \eqref{Psi} from 
Lemma \ref{lem:derivative of C}. 
For the higher order derivatives of 
amplitude functions, we have, for $|\mu|\ge 1$, 
\[
\left| \partial^{\mu}_{\xi}\alpha_{j,\pm}^k(\xi)\right| 
\le c, \quad 
\left| \partial^{\mu}_{\xi} \varepsilon_{j,\pm}^k(t;\xi)\right|  
\le c\, e^{\int^{|t|}_{0} 
(1+s)^{|\mu|}\Psi(s) \, ds}, \quad 
\quad |\xi| \ge 1,
\]
\[
\left| \partial^{\mu}_{\xi}\alpha_{j,\pm}^k(\xi)\right| 
\le c |\xi|^{-|\mu|}, \quad 
\left| \partial^{\mu}_{\xi} \varepsilon_{j,\pm}^k(t;\xi)\right| 
\le c\, e^{\int^{|t|}_{0} 
(1+s)^{|\mu|}\Psi(s) \, ds}|\xi|^{-|\mu|}, \quad 
0<|\xi|<1.
\]

If in addition to \eqref{L1}--\eqref{strict hyperbolicity2}, 
we further assume 
that $(1+\vert t \vert)^{|\mu|}a^\prime_{\nu,j}(t) \in 
L^1(\mathbb{R})$ for some $\mu$ with $|\mu|\ge 1$, and for all 
$\nu,j$, then the bound for each 
$\partial^{\mu}_{\xi} \varepsilon_{j,\pm}^k(t;\xi)$ is uniform in $t$.
\end{thm}

\noindent 
{\em Proof of Theorem \ref{thm:asymptotic}.} 
We must determine the precise 
growth order of $\alpha_{j,\pm}^k(\xi)$ and 
$\varepsilon_{j,\pm}^k(t;\xi)$ 
with respect to $\xi$. 

\begin{lem} \label{lem:lem3.1} 
Assume \eqref{L1}--\eqref{strict hyperbolicity2}. Then there exists a constant 
$c>0$ such that, for $j=1,\ldots,m$ and $k=0,\ldots,m-1$, 
\[
\left\vert \alpha^k_{j,\pm}(\xi)\right\vert \le c, 
\quad \left| \varepsilon^k_{j,\pm}(t;\xi)\right| \le c 
\int^{+\infty}_{\vert t \vert} \Psi(s) \, ds. 
\]
\end{lem}
\begin{proof} 
Let us go back to \eqref{Base}. Then we have, by using 
Lemma \ref{lem:lem2.2}, 
\begin{equation} \label{Prep}
\left|\pmb{\alpha}_{j,\pm}(\xi)+\pmb{\varepsilon}_{j,\pm}(t;\xi)\right|
\le \left|\Phi(t;\xi)^{-1}\mathcal{N}(t;\xi)\pmb{v}_j(t,\xi)\right|
\le c_0 
\end{equation}
for $j=1,\ldots,m$. 
Since $\pmb{\varepsilon}_{j,\pm}(t;\xi)$ decays to $\pmb{0}$, it follows that 
\[
\left|\pmb{\alpha}_{j,\pm}(\xi)\right|\le c_0. 
\]
On the other hand, \eqref{add} is equivalent to the following 
integral equation:
\[
\pmb{\varepsilon}_{j,\pm}(t;\xi)=i\int_{|t|}^{+\infty} 
C(s;\xi)\left( \pmb{\varepsilon}_{j,\pm}(s;\xi)
+\pmb{\alpha}_{j,\pm}(\xi)\right)\, ds.
\]
Thus combining this equation and \eqref{Prep}, we get 
\[
\left|\pmb{\varepsilon}_{j,\pm}(t;\xi)\right| \le c_0 \int_{|t|}^{+\infty} 
\|C(s;\xi)\| \, ds
\le c
\int^{+\infty}_{\vert t \vert} \Psi(s) \, ds. 
\]
The proof of Lemma \ref{lem:lem3.1} is finished. 
\end{proof} 

We need the estimates of higher order derivatives of amplitude functions. 
\begin{lem}\label{lem:useful lemma}
Assume \eqref{L1}--\eqref{strict hyperbolicity2}. Then we have, 
for $|\mu|\ge 1$, 
\[
\left| \partial^{\mu}_{\xi}\alpha^k_{j,\pm}(\xi)\right| 
\le c, \quad |\xi| \ge 1,
\]
\[
\left| \partial^{\mu}_{\xi}\alpha^k_{j,\pm}(\xi)\right| 
\le c |\xi|^{-|\mu|}, \quad 0<|\xi|<1,
\]
\begin{equation}
\left| \partial^{\mu}_{\xi} \varepsilon^k_{j,\pm}(t;\xi)\right| 
\le c\, e^{\int^{|t|}_{0} 
(1+s)^{|\mu|}\Psi(s) \, ds}, \quad |\xi|\ge1,
\label{small frequency8}
\end{equation}
\begin{equation}
\left| \partial^{\mu}_{\xi} \varepsilon^k_{j,\pm}(t;\xi)\right|
\le c\, e^{\int^{|t|}_{0} 
(1+s)^{\vert \mu \vert}\Psi(s) \, ds}
|\xi|^{-|\mu|}, \quad 0<|\xi|<1.
\label{small frequency9}
\end{equation}
 
In addition to \eqref{L1}--\eqref{strict hyperbolicity2}, if 
we assume that 
$(1+|t|)^{|\mu|} C(t;\xi) \in L^1(\mathbb{R})$ 
for some $\mu$ with $\vert \mu \vert \ge 1$, then 
\eqref{small frequency8}--\eqref{small frequency9} 
is uniform in $t$.
\end{lem} 
\begin{proof} 
Putting 
$$Q(t,\xi)=\left(\pmb{\alpha}_{0,\pm}(\xi)+\pmb{\varepsilon}_{0,\pm}(t;\xi),
\cdots,\pmb{\alpha}_{m-1,\pm}(\xi)+\pmb{\varepsilon}_{m-1,\pm}(t;\xi)
\right)
$$
we see that the matrix $Q(t,\xi)$ satisfies 
$$D_tQ(t,\xi)=C(t;\xi)Q(t,\xi)$$ with the initial data 
\[
Q(0,\xi)=\Phi(0;\xi)^{-1}\mathcal{N}(0;\xi)
(\pmb{v}_0(0;\xi),\cdots,\pmb{v}_{m-1}(0;\xi))=\mathcal{N}(0;\xi).
\]
Then it follows from the theory of ordinary differential 
equations that $Q(t,\xi)$ can be written by Picard series: 
\begin{multline}
Q(t,\xi)= \\
\p{I+i\int^{t}_0 C(\tau_1;\xi) \, d \tau_1
+i^2\int^{t}_0 C(\tau_1;\xi) \, 
d \tau_1\int^{\tau_1}_0 
C(\tau_2;\xi)\, d \tau_2+\cdots} \mathcal{N}(0;\xi).
\label{Picard series}
\end{multline}
We note from Lemma \ref{lem:derivative of C} that 
\begin{equation}
\left\|\partial_\xi^\mu C(t;\xi)\right\|
\le \begin{cases}
c(1+\vert t \vert)^{|\mu|}\Psi(t), 
& \text{$|\xi|\ge 1$,}\\ 
c(1+\vert t \vert)^{|\mu|}\Psi(t)
|\xi|^{-|\mu|}, & \text{$0<|\xi|<1$.}
\end{cases}
\label{C-Estimate}
\end{equation}
where 
$$\Psi(t)=\sum_{\underset{j\le m-1}{j+|\nu|=m}} 
|a^\prime_{\nu,j}(t)| 
\in L^1(\R).
$$ 
Differentiating \eqref{Picard series} with respect to $\xi$, 
we have, by using \eqref{C-Estimate}, 
\begin{equation}
\left| \partial^{\mu}_{\xi} 
\left(\pmb{\alpha}_{j,\pm}(\xi)+\pmb{\varepsilon}_{j,\pm}(t;\xi) \right) 
\right|
\le c\, e^{\int^{|t|}_0 
(1+s)^{\vert \mu \vert} \Psi(s) \, ds}
\label{another Picard}
\end{equation}
for all $t\in \mathbb{R}$, 
$\vert \xi \vert \ge 1$, $\vert \mu \vert \ge 1$ and 
$j=0,\ldots,m-1$, where we have used the following: 

\smallskip

\noindent 
{\bf Fact.} {\em 
Let $f(t) \in C(\mathbb{R})$. Then 
\[
e^{\int^{t}_{s} f(\tau) \, d \tau}
=1+\int^{t}_{s} f(\tau_1) \, d \tau_1
+\int^{t}_{s} f(\tau_1) \, 
d \tau_1 \int^{\tau_1}_{s} f(\tau_2) 
\, d \tau_2+\cdots.
\]
}

\smallskip 

Since 
$\partial_t \partial^{\mu}_{\xi} Q(t,\xi)=
\left(\partial_t \partial^{\mu}_{\xi} 
\pmb{\varepsilon}^{j}_{k,\pm}(t;\xi) \right)$, 
we combine $D_tQ=CQ$ and \eqref{another Picard} to deduce that 
$\partial_t\partial^{\mu}_{\xi}\pmb{\varepsilon}^{j}_\pm(0;\xi)$ exists 
and is uniformly bounded in $|\xi|\ge 1$, which ensures the 
existence of $\partial^{\mu}_{\xi}\pmb{\varepsilon}^{j}_\pm(0;\xi)$. 
Using again \eqref{another Picard} with $t=0$, 
we conclude 
\[
\left| \partial^{\mu}_{\xi}\pmb{\alpha}_{j,\pm}(\xi)\right| 
\le c, \quad |\xi|\ge1.
\]
In a similar way, we get the bound for 
$\partial^{\mu}_{\xi}\pmb{\alpha}_{j,\pm}(\xi)$ in low frequency 
part $0<|\xi|<1$. 
If we combine these estimates with 
\eqref{another Picard}, we have the bound for 
$\partial^{\mu}_{\xi}\pmb{\varepsilon}_{j,\pm}(t;\xi)$: 
\[
\left|\partial^{\mu}_{\xi}\pmb{\varepsilon}_{j,\pm}(t;\xi)\right|
\le \begin{cases}
c\, e^{\int^{|t|}_0 
(1+s)^{\vert \mu \vert} \Psi(s) \, ds},\quad & |\xi|\ge1,\\
c\, e^{\int^{|t|}_0 
(1+s)^{\vert \mu \vert} \Psi(s) \, ds}|\xi|^{-|\mu|}, \quad & 0<|\xi|<1.
\end{cases}
\]
The proof of Lemma \ref{lem:useful lemma} is complete. 
\end{proof}

\begin{proof}[Completion of the proof of Theorem \ref{thm:asymptotic}.] 
The estimates of the amplitude and error 
functions have been derived in 
Lemmas \ref{lem:lem3.1}--\ref{lem:useful lemma}. 
The proof of Theorem \ref{thm:asymptotic} is now finished.
\end{proof} 

\section{Estimates for oscillatory integrals; 
Proof of Theorem \ref{thm:Main theorem}}

The aim of this section is to establish time decay estimates
for $L^p$--$L^q$ norms of propagators for the Cauchy problem
\eqref{Equation}, which gives the proof of Theorem 
\ref{thm:Main theorem}. 
The analysis of high frequencies will give estimates
dependent on the geometry of the level sets of
characteristic roots of the equation.
For small frequencies estimates are
independent of the geometry of the level set and are given
by Proposition \ref{prop:low-freq} below. We recall that 
Theorem \ref{thm:asymptotic} assures in particular that the
solution to the Cauchy problem \eqref{Equation} is of the form
\begin{equation}
u(t,x)=\sum^{m-1}_{k=0}\sum^{m}_{j=1} \mathcal{F}^{-1} 
\left[ \left(\alpha^{j}_{k,\pm}(\xi)
+\varepsilon^{j}_{k,\pm}(t;\xi)\right)  
e^{i\vartheta_j(t;\xi)} \widehat{f}_k (\xi)\right](x), 
\quad t \gtrless 0.
\label{rep}
\end{equation}
The following Proposition \ref{prop:low-freq}
provides the time decay estimate
for small frequencies for each of the terms in this sum.
To simplify the notation, we formulate it in a more general 
form for general oscillatory integrals of the form
\[
T_t f(x)=\int_\Rn e^{i(x\cdot\xi+\vartheta(t;\xi))}
a(t,\xi) \widehat{f}(\xi) \, d\xi.
\] 
In the analysis of oscillatory integrals in the sum
\eqref{rep} we will actually  make 
time-dependent cut-offs and analyse separately
different ranges of frequencies.
We can obtain the following proposition
for small frequencies $|\xi|\leq t^{-1}$. Higher frequencies 
$|\xi|\geq t^{-1}$ will be analysed later. Thus, we
introduce a cut-off function of the form
$\psi((1+|t|)\xi)$ for some $\psi\in C_0^\infty(\Rn)$
such that $\psi(\xi)\equiv 1$ for 
$\vert \xi \vert \le \frac12$, 
and $0$ for  $\vert \xi \vert \ge 1$. 
We recall that we use the notation $\dot{L}^{p}_{\ka}(\Rn)$ for the
homogeneous Sobolev space $\dot{W}^{\kappa}_p(\Rn).$

\begin{prop}\label{prop:low-freq} 
Let $T_t$, $t\in \mathbb{R}$, be an operator defined by
\[
T_t f(x)=\int_\Rn e^{i(x\cdot\xi+\vartheta(t;\xi))}
\psi((1+|t|)\xi)a(t,\xi) \widehat{f}(\xi) \, d\xi,
\]
where $\vartheta(t;\xi)$ is real valued, 
positively homogeneous of order one in 
$\xi$. Assume that the amplitude $a(t,\xi)$ satisfies
\begin{equation}\label{bdd of a}
\vert \xi \vert^\ka | a(t,\xi)|\leq C 
\end{equation}
for some $\ka\in\R$, and for all $t\in\R$ and all 
$\xi\in\supp \psi((1+|t|)\xi)$. Let 
$1\le p\leq 2\leq q \le +\infty$ be such that
$\frac{1}{p}+\frac{1}{q}=1$.
Then for $t \in \mathbb{R}$ we have the estimate
\begin{equation}\label{small Lp-Lq}
 \Vert T_t f\Vert_{L^q(\Rn)} \leq 
 C(1+\vert t\vert)^{-n\left(\frac{1}{p}-\frac{1}{q}\right)}
 \Vert f \Vert_{\dot{L}^{p}_{-\ka}(\Rn)}, 
\end{equation}
where constant $C$ depends on $n,p,q$ and the
norm $\left\Vert \vert \xi \vert^\ka a \right\Vert_{L^\infty}$.
\end{prop}
\begin{proof} We easily obtain
$$\|T_tf\|_{L^2(\Rn)}\le C\|f\|_{L^2_{-\kappa}(\Rn)}
$$ 
by the Plancherel identity.
Thus \eqref{small Lp-Lq}
would follow by analytic interpolation 
from an estimate:
\begin{equation}\label{EQ:Main EST}
\|T_tf\|_{L^\infty(\Rn)}
\le C(1+|t|)^{-n}\|f\|_{\dot{L}^1_{-\kappa}(\Rn)}. 
\end{equation} 
In fact, since $|a(t,\xi)| |\widehat{f}(\xi)|
\le C\||\xi|^\kappa a\|_{L^\infty} 
\|f\|_{L^1_{-\kappa}(\Rn)}$, we can estimate 
\begin{align*}
& \|T_tf\|_{L^\infty(\Rn)}
\le \int_{|\xi|\le (1+|t|)^{-1}}
|a(t,\xi)| |\widehat{f}(\xi)|\, d\xi\\
\le& C\left(\int_{|\xi|\le (1+|t|)^{-1}}d\xi\right)\|f\|_{L^1_{-\kappa}(\Rn)}
\le C(1+|t|)^{-n} \|f\|_{\dot{L}^1_{-\kappa}(\Rn)},
\end{align*}
for all $t\in \R$. 
This proves \eqref{EQ:Main EST}. 
The proof of Proposition \ref{prop:low-freq} 
is complete. 
\end{proof}

We now turn to the analysis of larger frequencies.
The following proposition provides the necessary 
background to obtain the time decay estimate
for large frequencies for each of the terms in the sum 
\eqref{rep}. In fact, it will be used for
frequencies $|\xi|\geq 1$ but will be formulated here in
a slightly more general form. The relation with the sum
\eqref{rep} and a refinement for low frequencies 
will be also made in Proposition
\ref{COR:oscillatory-convex}.

\begin{prop} \label{PROP:oscillatory-convex}
Let $T_t$, $t>0$, be an operator defined by
\[
T_t f(x)=\int_\Rn e^{i(x\cdot\xi+\vartheta(t;\xi))}
a(t,\xi) \widehat{f}(\xi) \, d\xi,
\]
where $\vartheta(t;\xi)$ is real valued, continuous in $t$, 
smooth in
$\xi\in\Rn\backslash 0$, homogeneous of order one in $\xi$.
Assume that the set
$$\Sigma_\varphi=\{\xi\in\Rn\backslash 0:\va(\xi)=1\}$$
is strictly convex and let 
$\gamma=\gamma(\Sigma_\varphi)$ be
the convex Sugimoto index of $\Sigma_\varphi$,
as defined in {\rm\eqref{eq:convex-sugimoto}}. 
Suppose that 
$$|\vartheta (t;\xi)|\leq 
C(1+t)|\xi|\quad \text{for all}
\quad t>0,\,\xi\in\Rn,$$
and that there is some $\varphi\in C^\infty(\Rn\backslash 0)$,
$\va>0$, 
such that 
\begin{equation}\label{EQ:conv-theta}
 t^{-1}\partial_\xi^\alpha\vartheta(t;\xi)\to 
\partial_\xi^\alpha\varphi(\xi)
 \textrm{ as } t\to\infty,
 \textrm{ for all }\; \xi\in\Rn\backslash 0,\;
|\alpha|\leq\gamma.
\end{equation} 
Assume also that the amplitude $a(t,\xi)$ satisfies
\begin{equation}\label{EQ:assAMP}
|\partial_\xi^\alpha a(t,\xi)|\leq C_\alpha
\quad\text{for all}
\quad |\alpha|\leq [(n-1)/\gamma]+1.
\end{equation}
Let $1< p\leq 2\leq q< \infty$ be such that
$\frac{1}{p}+\frac{1}{q}=1$.
Then for $t>0$ we have the estimate
\begin{equation}\label{EQ:convex-estimate}
 \Vert T_t f \Vert_{L^q(\Rn)} 
 \leq C t^{-\frac{n-1}{\gamma}\p{\frac{1}{p}-\frac{1}{q}}}
 \Vert f \Vert_{{L}^p_{N_p}(\Rn)},
\end{equation}
where $N_p=\left(n-\frac{n-1}{\gamma}+
\left[ \frac{n-1}{\gamma} \right]+1\right)
\left(\frac{1}{p}-\frac{1}{q}\right)$. 
\end{prop}
The number of derivatives $N_p=\left(n-\frac{n-1}{\gamma}+
\left[ \frac{n-1}{\gamma} \right]+1\right)
\left(\frac{1}{p}-\frac{1}{q}\right)$ required for the
estimate \eqref{EQ:convex-estimate} is determined by the fact
that the amplitude $a(t,\xi)$ in \eqref{EQ:assAMP}
is in the symbol class $S^0_{0,0}$ rather than the usual
$S^0_{1,0}$. In fact, if $a(t,\xi)$ satisfies inequalities
\begin{equation}\label{EQ:assAMP2}
|\partial_\xi^\alpha a(t,\xi)|\leq C_\alpha \jp{\xi}^{-|\alpha|}
\quad\text{for all}
\quad |\alpha|\leq \left[\frac{n-1}{\gamma}\right]+1,
\end{equation}
then we can take the Sobolev index
$N_p=\big(n-\frac{n-1}{\gamma}
\big)\big(\frac{1}{p}-\frac{1}{q}\big)$ for 
the estimate \eqref{EQ:convex-estimate} to hold.
However, the method of asymptotic integration and the
statement of Theorem \ref{thm:asymptotic} forces us to
assume \eqref{EQ:assAMP} rather than \eqref{EQ:assAMP2}.

Let us now discuss other assumptions we make in this proposition
from the point of view of the original Cauchy problem
\eqref{Equation}.
We recall from \eqref{Equation-lim} that
functions $\varphi^\pm_k(\xi)$ are characteristics of the 
limiting strictly hyperbolic operator
\begin{equation}
L^\pm(D_t,D_x)u 
\equiv D^{m}_t u+\sum_{\underset{j \le m-1}
{\vert \nu \vert+j=m}} 
a^\pm_{\nu,j} D^{\nu}_x D^{j}_t u,
\label{Equation-lim1}
\end{equation}
and their geometric properties are responsible for the time decay
of solutions to the Cauchy problems for both operators
$L(t,D_t,D_x)$ and $L^\pm(D_t,D_x)$. 
The fact that $\va_k^\pm(\xi)$ are characteristics of
\eqref{Equation-lim1}, implies that they are 
real analytic for $\xi\not=0$ and that we have
the following statement, which was established for operators
with constant coefficients by Sugimoto \cite{sugi94}.
\begin{prop}
Let $\va_k(\xi)$, $k=1,\ldots,m$, be characteristics of
operator \eqref{Equation-lim1}, ordered by
$\va_1(\xi)>\va_2(\xi)>\cdots>\va_m(\xi)$ for $\xi\not=0$.
Suppose that all the Hessians $\va_k^{\prime\prime}(\xi)$ are
semi-definite for $\xi\not=0$. Then there exists a polynomial
$\alpha(\xi)$ of order one such that
$\va_{m/2}(\xi)>\alpha(\xi)>\va_{m/2+1}$ 
{\rm (}if $m$ is even{\rm )} or
$\alpha(\xi)=\va_{(m+1)/2}(\xi)$ {\rm (}if $m$ is odd{\rm )}. Moreover,
the hypersurfaces $\Sigma_k=\{\xi\in\Rn;\, 
\widetilde{\va}_k=\pm 1\}$
with $\widetilde{\va}_k(\xi)=\va_k(\xi)-\alpha(\xi)$ 
$(k\not=(m+1)/2)$ are convex and $\gamma(\Sigma_k)\leq 2[m/2].$
\label{PROP:limitingphases}
\end{prop}
In particular, in our arguments we can replace $\va_k$ by
$\widetilde{\va}_k$ since the addition of a linear function
does not change the decay rate nor the index 
$\gamma(\Sigma_{\va_k}).$
This also ensures that the limiting phase $\va$ in
Proposition \ref{PROP:oscillatory-convex} may be taken
to be strictly positive. Indeed, it can be taken to
be nonzero, and if it is strictly negative we simply
replace $\va$ by $-\va$.
Moreover, the
assumption that $\Sigma_\va$ is strictly convex 
in Proposition \ref{PROP:oscillatory-convex}
can be replaced by the assumption that it is only convex. 
Indeed, 
since $\va(\xi)$ is a characteristic root of
\eqref{Equation-lim1}, it is real analytic for $\xi\not=0$.
Then, the convexity, the real analyticity and the compactness
imply that it is actually strictly convex. In particular,
it also implies that $\gamma(\Sigma_\va)$ is {\em finite and even}.

In the case when the level set $\Sigma_\varphi$ in Proposition
\ref{PROP:oscillatory-convex} is not convex, we get a weaker
result based on the one-dimensional
van der Corput lemma. In this case we
use the non-convex Sugimoto index of $\Sigma_\varphi$ which
was defined in \eqref{eq:nonconvex-Sugimoto} in the introduction.

\begin{prop} \label{PROP:oscillatory-nonconvex}
Let $T_t$, $t>0$, be an operator defined by
\[
T_t f(x)=\int_\Rn e^{i(x\cdot\xi+\vartheta(t;\xi))}
a(t,\xi) \widehat{f}(\xi) \, d\xi,
\]
where $\vartheta(t;\xi)$ is real valued, continuous in $t$, 
smooth in
$\xi\in\Rn\backslash 0$, 
homogeneous of order one in $\xi$.
Let $\gamma_0=\gamma_0(\Sigma_\va)$ be the
non-convex Sugimoto index of the level surface
$\Sigma_\varphi=\{\xi\in\Rn\backslash 0:\va(\xi)=1\}$.
Suppose that 
$$|\vartheta (t;\xi)|\leq C(1+t)|\xi|\quad \text{for all}
\quad t>0,\,\xi\in\Rn,$$
and that there is some $\varphi\in C^\infty(\Rn\backslash 0)$,
$\va>0$,
such that 
$$t^{-1}\partial_\xi^\alpha\vartheta(t;\xi)\to 
\partial_\xi^\alpha\varphi(\xi)
\quad \text{as}\quad t\to\infty,
\quad \text{for all}\quad \xi\in\Rn,\; |\alpha|\leq\gamma_0.
$$ 
Assume also that the amplitude $a(t,\xi)$ satisfies
\begin{equation}\label{EQ:assAMP-nc}
|\partial_\xi^\alpha a(t,\xi)|\leq C_\alpha
\quad\text{for all}
\quad |\alpha|\leq 1.
\end{equation}
Let $1< p\leq 2\leq q<+\infty$ be such that
$\frac{1}{p}+\frac{1}{q}=1$.
Then for $t>0$ we have the estimate
\begin{equation}\label{EQ:convex-estimate-nc}
 \Vert T_t f \Vert_{L^q(\Rn)} \leq C t^{-\frac{1}{\gamma_0}\p{\frac{1}{p}-
 \frac{1}{q}}}
 \Vert f \Vert_{{L}^p_{N_p}(\Rn)},
\end{equation}
where $N_p=\left(n-\frac{1}{\gamma_0}+1\right)
\left(\frac{1}{p}-\frac{1}{q}\right)$. 
\end{prop}
We will first prove Proposition {\rm \ref{PROP:oscillatory-convex}}
and then indicate the changes necessary for the proof of
Proposition \ref{PROP:oscillatory-nonconvex}.
\begin{proof}[Proof of Proposition {\rm \ref{PROP:oscillatory-convex}}]
First we observe that $\Vert T_t f \Vert_{L^2}\leq C \Vert f \Vert_{L^2}$
by the Plancherel identity. 
In order to simplify the proof somewhat, we will absorb the
Sobolev index $N_p$ into the amplitude $a(t,\xi)$, so that we will
estimate $||T_tf||_{L^q}$ in terms of $||f||_{L^p}$, and not
in terms of $||f||_{L^p_{N_p}}$. Thus, instead of 
\eqref{EQ:assAMP}, from now on
we will assume that
the amplitude $a(t,\xi)$ satisfies
$$|\partial_\xi^\alpha a(t,\xi)|\leq C_\alpha \jp{\xi}^{-k}$$ for
all $|\alpha|\leq [(n-1)/\gamma]+1$ and $k=N_p$, so that
estimate
\eqref{EQ:convex-estimate} would follow by interpolation
from the estimate
\begin{equation}\label{EQ:est2}
\Vert T_t f \Vert_{L^\infty(\Rn)} \leq C t^{-\frac{n-1}{\gamma}}
 \Vert f \Vert_{L^1(\Rn)},
\end{equation}
and where we take $k=N_1=n-\frac{n-1}{\gamma}+
\left[ \frac{n-1}{\gamma} \right]+1$.
Note that since we assume that the amplitude is bounded for small
frequencies, we can work with standard Sobolev spaces here.
By using Besov spaces, we can microlocalise the desired estimate
to discs in the frequency space. Indeed, let 
$\{\Phi_j\}_{j=0}^\infty$ be the Littlewood-Paley partition of
unity, and let 
$$
\Vert u \Vert_{B^s_{p,q}}=\p{ \sum_{j=0}^\infty \p{2^{js} \Vert \Fcal^{-1}
\Phi_j(\xi)\Fcal u \Vert_{L^q(\Rn)}}^q}^{1/p}
$$
be the norm of the Besov space $B^s_{p,q}$. Then, because of the
continuous embeddings $L^p\subset B^0_{p,2}$ for $1<p\leq 2$,
and $B^0_{q,2}\subset L^{q}$ for $2\leq q<+\infty$ 
(see \cite{BL}), it is sufficient to prove the uniform estimate
for the operators with amplitudes $a(t,\xi)\Phi_j(\xi)$.
Let us denote 
$$\widetilde{\vartheta}(t,\xi)=t^{-1}\vartheta(t,\xi),$$
so that by the assumption we have
$\widetilde{\vartheta}(t;\xi)\to\varphi(\xi)$ as $t\to\infty$.
Now, writing 
$$\Phi_j(\xi)=\Phi_j(\xi)\Psi\p{\frac{
\widetilde\vartheta(t;\xi)}{2^j}}$$ with
some function $\Psi\in C_0^\infty(0,\infty)$, we may prove
the uniform estimate for operators with amplitudes
$a(t,\xi)\Psi\p{\frac{
\widetilde\vartheta(t;\xi)}{2^j}}.$
Such choice of $\Psi$ is possible due to our assumption that
$\widetilde{\vartheta}(t;\xi)\to\varphi(\xi)$ as $t\to\infty$,
and we restrict the analysis for large enough $t$.
Let 
\begin{equation}\label{EQ:Itx}
I(t,x)=\int_\Rn e^{i(x\cdot\xi+\vartheta(t;\xi))}
a(t,\xi)\Psi\p{\frac{\widetilde\vartheta(t;\xi)}{2^j}}\, d\xi
\end{equation}
be the kernel of the corresponding operator.
Since we easily have the $L^2$--$L^2$ estimate by the 
Plancherel identity, by analytic interpolation
we only need to prove the $L^1$--$L^\infty$ case of \eqref{EQ:est2}.
In turn, this follows from the
estimate $|I(t,x)|\leq C t^{-\frac{n-1}{\gamma}},$
with constant $C$ independent of $j$.

Let $\kappa\in C_0^\infty(\Rn)$ be supported in 
a ball with some radius $r>0$ centred at the origin. 
We split the integral in
\begin{align*}
I(t,x) & =  I_1 (t,x)+I_2 (t,x) \\ & =
\int_\Rn e^{i(x\cdot\xi+\vartheta(t;\xi))}
a(t,\xi) \kappa\left(t^{-1}x+t^{-1}\nabla_\xi\vartheta(t;\xi)\right)
\Psi\p{\frac{\widetilde\vartheta(t;\xi)}{2^j}} d\xi
\\
& + 
\int_\Rn e^{i(x\cdot\xi+\vartheta(t;\xi))}
a(t,\xi) (1-\kappa)
\left(t^{-1}x+t^{-1}\nabla_\xi\vartheta(t;\xi)\right)
\Psi\p{\frac{\widetilde\vartheta(t;\xi)}{2^j}} d\xi.
\end{align*}
We can easily see that 
$|I_2(t,x)|\leq C t^{-\frac{n-1}{\gamma}}$. In fact,
we can show
$|I_2(t,x)|\leq C t^{-l}$ 
for $l=[(n-1)/\gamma]+1$ and then the required estimate
simply follows since $l>(n-1)/\gamma$. 
Indeed, on the support of $1-\kappa$,
we have $|x+\nabla_\xi\vartheta (t;\xi)|\geq rt>0$. Thus,
integrating by parts with operator
$P=\frac{x+\nabla_\xi \vartheta (t;\xi)}
{i|x+\nabla_\xi\vartheta (t;\xi)|^2}
\cdot \nabla_\xi$, we get
\begin{multline}\label{EQ:I2tx}
I_2(t,x)= \\ \int_\Rn 
e^{i(x\cdot\xi+\vartheta(t;\xi))}
(P^*)^l \left[ a(t,\xi) (1-\kappa)
\left(t^{-1}x+t^{-1}\nabla_\xi\vartheta(t;\xi)\right)
\Psi\p{\frac{\widetilde\vartheta(t;\xi)}{2^j}}\right]  d\xi.
\end{multline}
Using the fact that $|\partial_\xi^\alpha\vartheta(t;\xi)|
\leq C(1+t)|\xi|^{1-|\alpha|}$, 
we readily observe from \eqref{EQ:I2tx} that 
the required estimate $|I_2(t,x)|\leq C t^{-l}$ holds. Here we
also used the condition \eqref{EQ:assAMP} which assures that
we can perform the integration by parts $[(n-1)/\gamma]+1$
times.

Now we will turn to estimating $I_1(t,x)$. Recall
that $\widetilde{\vartheta}(t;\xi)=t^{-1}\vartheta(t;\xi)$ and
$\widetilde{\vartheta}(t;\xi)\to\varphi(\xi)$ as $t\to\infty$.
Let us denote 
$$\Sigma^t=\{\xi\in\Rn: 
\widetilde{\vartheta}(t;\xi)=1\}.$$ 
It can be readily checked
that $\gamma(\Sigma^t)\to \gamma(\Sigma_\varphi)=\gamma$ as 
$t\to\infty$. So we can restrict our attention to $t$ large
enough for which we have $\gamma(\Sigma^t)=\gamma$.
By rotation, we can always microlocalise in some narrow cone
around $e_n=(0,\ldots,0,1)$ and in this cone we can parameterise 
$$\Sigma^t=\{(y,h_t(y)): y\in U\}$$ 
for some open $U\subset\R^{n-1}$.
In other words, we have 
$\widetilde{\vartheta}(t;y,h_t(y))=1$, and it follows that
$h_t$ is smooth and
$\nabla h_t:U\to \nabla h_t(U)\subset\R^{n-1}$ is a
homeomorphism. The function $h_t$ is concave if $\Sigma^t$ is convex.
We claim that 
\begin{equation}\label{EQ:ht}
|\partial_y^\alpha h_t(y)|\leq C_\alpha,\quad
\textrm{for all}\quad y\in U\quad \textrm{and large enough}
\quad t.
\end{equation}
Indeed, let us look at $|\alpha|=1$ first.
From $\widetilde{\vartheta}(t;y,h_t(y))=1$ we get that
$$\nabla_y \widetilde{\vartheta}+\partial_{\xi_n}
\widetilde{\vartheta}\cdot \nabla h_t(y)=0.$$
From homogeneity we have 
$|\nabla_\xi\widetilde{\vartheta}|\leq C$,
so also $|\nabla_y \widetilde{\vartheta}|\leq C$.
By Euler's identity we have
\begin{equation}\label{EQ:phi-non-zero}
\partial_{\xi_n} \widetilde{\vartheta} (t;e_n)=
\widetilde{\vartheta}(t;e_n)\to \varphi(e_n)>0
\textrm{ as }t\to\infty,
\end{equation}  
so we have $|\partial_{\xi_n} \widetilde{\vartheta}|\geq c>0$
since we are in a narrow cone around $e_n$. From this it
follows that $|\nabla_y h_t(y)|\leq C$ for all $y\in U$
and $t$ large enough. A similar argument proves the boundedness
of higher order derivatives in \eqref{EQ:ht}.

Now, let us turn to analyse the structure of the sets $\Sigma^t$.
We have the Gauss map
$$
\nu:\Sigma^t\ni\zeta \mapsto \frac{
\nabla_\zeta\widetilde{\vartheta}(t;\zeta)}
{|\nabla_\zeta\widetilde{\vartheta}(t;\zeta)|}\in\Snm,
$$
and for $x=(x^\prime,x_n)\in\R^{n-1}\times\R$
near the point $-\nabla_\zeta\widetilde{\vartheta}(t;e_n)$
we define $z_t\in U$ by
$$(z_t,h_t(z_t))=\nu^{-1}(-x/|x|).$$ 
Then $(-\nabla_y h_t(y),1)$ is normal to $\Sigma^t$ at
$(y,h_t(y))$, so we get
$$
-\frac{x}{|x|}=\frac{(-\nabla_y h_t(z_t),1)}
{|(-\nabla_y h_t(z_t),1)|} \quad \textrm{and} \quad
\frac{x^\prime}{x_n}=-\nabla_y h_t(z_t).
$$
Making change of variables $\xi=(\widetilde\lambda y,
\widetilde\lambda h_t(y))$
and using $\widetilde{\vartheta}(t;\xi)=\widetilde\lambda$, we get
\begin{eqnarray}
I_1(t,x) & =&
\int_0^\infty \int_U e^{i\widetilde\lambda(x^\prime\cdot y+x_n h_t(y)+t)}
a(t,\widetilde\lambda y,\widetilde\lambda h_t(y))
\Psi\p{\frac{\widetilde\lambda}{2^j}}\kappa_0(t,x,y) 
\abs{\frac{d\xi}{d(\widetilde\lambda,y)}} \, dy d\widetilde\lambda 
\nonumber \\
& = &
 \int_0^\infty \int_U 
e^{i\widetilde\lambda(-x_n\nabla_y h_t(z_t)\cdot y+x_n h_t(y)+t)}
\left[\widetilde\lambda^l a(t,\widetilde\lambda y,
\widetilde\lambda h_t(y))\right]
\widetilde\lambda^{n-1-l} \times \nonumber \\
& &\quad \times
\Psi\p{\frac{\widetilde\lambda}{2^j}}\kappa_0(t,x,y) \chi(t,y) 
\, dy d\widetilde\lambda
 \label{EQ:I1tx} \\
& = & 
\int_0^\infty \int_U 
e^{i\lambda(-\nabla_y h_t(z_t)\cdot y+ h_t(y)+t x_n^{-1})}
\widetilde{a}(t,x_n,\lambda y,\lambda h_t(y)) 
\lambda^{n-1-l}  \times \nonumber \\
& & \quad \times \Psi\p{\frac{\lambda}{2^j x_n }} x_n^{-n+1+l-1}
\kappa_0(t,x,y) \chi(t,y) \, dy d\lambda,\nonumber 
\end{eqnarray}
where $\kappa_0(t,x,y)=
\kappa\left(t^{-1}x+\nabla_\xi\widetilde{\vartheta}(t;y,h_t(y))\right)$,
and 
$$\widetilde{a}(t,x_n,\lambda y,\lambda h_t(y))=
(x_n^{-1}\lambda)^l a
\left(t,x_n^{-1} \lambda y,x_n^{-1} \lambda h_t(y)\right),$$
and
where we made a change $\widetilde\lambda=x_n^{-1}\lambda$ in the
last equality. Here also we used
$\abs{\frac{d\xi}{d(\widetilde\lambda,y)}}=
\widetilde\lambda^{n-1} \chi(t,y)$,
where $\chi(t,y)$ and all of its derivatives with respect to $y$
are bounded because of \eqref{EQ:ht}.

If we choose $r$ in the definition of the cut-off function
$\kappa$ sufficiently small, then on its support we have 
$|x|\approx |x_n| \approx t$, and we can estimate
\begin{equation}\label{EQ:IJ}
\begin{aligned}
|I_1(t,x)| & \leq C t^{-n+l}\int_0^\infty 
\abs{J(\lambda,z_t)
\Psi\p{\frac{\lambda}{2^j t }}
\lambda^{n-1-l}} \, d\lambda \\ 
& =
C t^{-n+l}2^{j(n-l)}\int_0^\infty 
\abs{J(2^j \lambda,z_t)
\Psi\p{\frac{\lambda}{t }}
\lambda^{n-1-l}} \, d\lambda
\end{aligned}
\end{equation}
with
$$
J(\lambda,z_t)=\int_U 
e^{i\lambda(-\nabla_y h_t(z_t)\cdot y+ h_t(y)+t x_n^{-1})}
\widetilde{a}(t,x_n,\lambda y,\lambda h_t(y)) 
\kappa_0(t,x,y) \chi(t,y)\, dy.
$$
We will show that
\begin{equation}\label{EQ:estJ}
\abs{J(\lambda,z_t)}\leq C(1+\lambda)^{-\frac{n-1}{\gamma}},
\quad \lambda>0.
\end{equation}
Then, if we take $l=n-\frac{n-1}{\gamma}$, and use
\eqref{EQ:IJ} and \eqref{EQ:estJ}, we get
\begin{equation}\label{EQ:I1tx-est2}
\begin{aligned}
|I_1(t,x)| & \leq C t^{-\frac{n-1}{\gamma}} 
2^{j\frac{n-1}{\gamma}}
\int_0^\infty
(2^j\lambda)^{-\frac{n-1}{\gamma}}
\Psi\p{\frac{\lambda}{t}}\lambda^{\frac{n-1}{\gamma}-1} \, d\lambda
\\
& = C t^{-\frac{n-1}{\gamma}} 
\int_0^\infty
\lambda^{-1}
\Psi\p{\frac{\lambda}{t}} \, d\lambda 
=  C t^{-\frac{n-1}{\gamma}} 
\int_0^\infty
\lambda^{-1}
\Psi\p{\lambda}\, d\lambda \\
& \leq C t^{-\frac{n-1}{\gamma}},
\end{aligned}
\end{equation}
which is the desired estimate
for $I_1(t,x)$.

Let us now prove \eqref{EQ:estJ}. 
It will, in turn, follow from Theorem \ref{THM:oscintthm} below.
First of all we note that since we assumed that
$|\partial_\xi^\alpha a(t,\xi)|\leq C_\alpha \jp{\xi}^{-N_1}$ for
all $|\alpha|\leq [(n-1)/\gamma]+1$ and 
$N_1=n-\frac{n-1}{\gamma}+
\left[ \frac{n-1}{\gamma} \right]+1$, we get that
\begin{equation}\label{eq:a-tilde}
 |\partial^\alpha_y \widetilde{a}|\leq C \quad
 \textrm{for all}\quad |\alpha|\leq 
 \left[ (n-1)/\gamma \right]+1.
\end{equation}

To write $J(\lambda,z_t)$ in a suitable form,
we change to polar coordinates
$(\rho,\omega)$ with $y=\rho\omega+z_t$, so that
\begin{equation}\label{EQ:gather1}
J(\lambda,z_t)=\int_{{\mathbb S}^{n-2}} \int^{\infty}_{0}
e^{i\lambda F(\rho,z_t,\omega)} \beta(\rho,z_t,\omega) 
\rho^{n-2} \, d\rho d\omega,
\end{equation}
with
\begin{gather}\label{osc-phase}
F(\rho,z_t,\omega)= h_t(\rho\omega+z_t)-h_t(z_t)-
\rho \nabla_y h_t(z_t)\cdot\omega, \\
\label{EQ:gather3}
\beta(\rho,z_t,\omega)=
\widetilde{a}\left(t,x_n,\lambda (\rho\omega+z_t),\lambda h_t
(\rho\omega+z_t)\right) 
\kappa_0(t,x,\rho\omega+z_t) \chi(t,\rho\omega+z_t),
\end{gather}
where we can assume in addition that $\chi=0$ unless
$\rho\omega+z_t\in U$, so both $\rho$ and $\omega$ vary over
bounded sets.

Now we can apply the following result,
which has appeared in \cite{R-vdC}
for more general complex valued phases $\Phi$, 
thus including the real-valued case of the 
phase function $F$ in \eqref{osc-phase}.
The estimate \eqref{EQ:estJ} follows from the following
theorem with $N=n-1$. 

\begin{thm}[\cite{R-vdC}]\label{THM:oscintthm}
Consider the oscillatory integral
\begin{equation*}\label{EQ:genoscint}
I(\la,\nu)=\int_{\R^{N}}
e^{i\la\Phi(x,\nu)}a(x,\nu)\chi(x)\,dx\,,
\end{equation*}
where $N\geq 1$\textup{,} 
and $\nu$ is a parameter. 
Let $\ga\ge2$ be an integer.
Assume that
\begin{itemize}
\item[(A1)]\label{HYP:mainoscintgbdd} there exists a 
sufficiently small $\de>0$ such that 
$\chi\in C^\infty_0(B_{\delta/2}(0))$, where
$B_{\delta/2}(0)$ is the ball with radius 
${\delta/2}$ around $0$\textup{;}
\item[(A2)]\label{HYP:mainoscintImPhipos} $\Phi(x,\Ivar)$ is a
complex valued function such that $\Im \Phi(x,\Ivar)\ge0$ for all
$x\in \supp\chi$ and all parameters $\Ivar$\textup{;}
\item[(A3)]\label{HYP:mainoscintFconvexfn} for some fixed
$z\in\supp\chi$,  the function
\begin{equation*}
F(\rho,\om,\Ivar):=\Phi(z+\rho\om,\Ivar), \; |\om|=1,
\end{equation*}
satisfies the following conditions. Assume that
for each $\mu=(\omega,\nu)$, function $F(\cdot,\mu)$ is of class
$C^{\gamma+1}$ on $\supp\chi$, and
let us write its $\gamma^{\rm th}$ order
Taylor expansion in $\rho$ at $0$ as
\begin{equation*}\label{EQ:Fformwithremainder}
F(\rho,\mu)=\sum_{j=0}^\gamma a_j(\mu)\rho^j + 
R_{\gamma+1}(\rho,\mu)\,,
\end{equation*}
where $R_{\gamma+1}$
is the remainder term. Assume that 
we have
\begin{itemize}
\item[(F1)]\label{ITEM:AssumpF1} $a_0(\mu)=a_1(\mu)=0$ for all
$\mu$\textup{;}
\item[(F2)]\label{ITEM:AssumpF2} there exists a constant $C>0$ such that
$\sum_{j=2}^\ga\abs{a_j(\mu)}\ge C$ for all
$\mu$\textup{;}
\item[(F3)]\label{ITEM:AssumpF3} for each
$\mu$\textup{,} $\abs{\pa_\rho F(\rho,\mu)}$ is
increasing in $\rho$ for $0<\rho<\de$\textup{;}
\item[(F4)]\label{ITEM:AssumpF4} for each $k\leq \gamma+1$,
$\pa_\rho^kF(\rho,\mu)$ is bounded uniformly in
$0<\rho<\de$ and $\mu$\textup{;}
\end{itemize}
\item[(A4)]\label{HYP:mainoscintAderivsbdd} for each multi-index $\al$
of length $\abs{\al}\le \big[\frac{N}{\ga}\big]+1$\textup{,} there
exists a constant $C_\al>0$ such that 
$\abs{\pa_x^\al a(x,\Ivar)}\le
C_\al$ for all $x\in \supp\chi$ and all parameters $\Ivar$.
\end{itemize}
Then there exists a constant $C=C_{N,\ga}>0$ such that
\begin{equation}\label{EQ:oscintbound}
\abs{I(\la,\Ivar)}\le C(1+\la)^{-\frac{N}{\ga}}\quad\text{for all
}\; \la\in[0,\infty)
\textrm{ and all parameters } \Ivar.
\end{equation}
\end{thm}
We refer to \cite{R-vdC} and
\cite{RS2} for details. 
Now, the function $F$ in \eqref{osc-phase} satisfies
condition (A3) of Theorem \ref{THM:oscintthm}
because of the definition of the
convex Sugimoto index $\gamma$ and because $h_t$ is concave.
Since $\partial_y^\alpha h_t$, $|\alpha|\leq\gamma$,
can be expressed via $\partial_\xi^\alpha \varphi$,
$|\alpha|\leq\gamma$, and since we have
\eqref{EQ:phi-non-zero}, it also follows from
\eqref{EQ:convphi} and \eqref{EQ:conv-theta} that 
function $F$ satisfies property (F2) of
Theorem \ref{THM:oscintthm}.
The proof of Proposition \ref{PROP:oscillatory-convex} 
is now complete.
\end{proof}

Let us now show that we can actually also insert the cut-off
$1-\psi((1+\vert t \vert)\xi)$ in Proposition 
\ref{PROP:oscillatory-convex}
which is necessary for the analysis of the representation
\eqref{rep}. Here $\psi\in C_0^\infty(\Rn)$
such that $\psi(\xi)\equiv 1$ for $\vert \xi \vert \le \frac12$, 
and $0$ for  $\vert \xi \vert \ge 1$. The case of high frequencies
$|\xi|\geq 1$ (for solutions) is covered by Proposition
\ref{PROP:oscillatory-convex}, and the proof about the
insertion of $1-\psi((1+\vert t \vert)\xi)$ is similar
to the proof of the following Proposition 
\ref{COR:oscillatory-convex}.
So we now restrict to
$t^{-1}<|\xi|\leq 1$, since the case $|\xi|<t^{-1}$ was covered
in Proposition \ref{prop:low-freq}.

\begin{prop} \label{COR:oscillatory-convex}
Let $T_t$, $t \ne 0$, be an operator defined by
\begin{equation}\label{EQ:oscillatory-cor}
T_t f(x)=\int_\Rn e^{i(x\cdot\xi+\vartheta(t;\xi))}
\br{1-\psi((1+|t|)\xi)}
a(t,\xi) \widehat{f}(\xi) \, d\xi,
\end{equation}
where $\vartheta(t;\xi)$, and 
$\gamma$ are as in Proposition
\ref{PROP:oscillatory-convex}. Assume that the amplitude
$a(t,\xi)$ satisfies $a(t,\xi)=0$ for all $|\xi|\geq 1$ and 
that 
\[
|\partial_\xi^\alpha a(t,\xi)|\leq C_\alpha |\xi|^{-|\alpha|}
\quad\text{for all}
\quad |\alpha|\leq [(n-1)/\gamma]+1.
\]
Let $1< p\leq 2\leq q< \infty$ be such that
$\frac{1}{p}+\frac{1}{q}=1$.
Then for $t\ne 0$ we have the estimate
\begin{equation}\label{EQ:convex-estimate-cor}
 \Vert T_t f \Vert_{L^q(\Rn)} \leq 
 C \vert t\vert^{-\frac{n-1}{\gamma}\p{\frac{1}{p}-
 \frac{1}{q}}}
 \Vert f \Vert_{\dot{L}^p_{N_p}(\Rn)},
\end{equation}
where $N_p=\left(n-\frac{n-1}{\gamma}
+\left[ \frac{n-1}{\gamma} \right]+1
\right)
\left(\frac{1}{p}-\frac{1}{q}\right)$. 
\end{prop}
\begin{proof}
The proof of this proposition is almost the same as the proof of
Proposition \ref{PROP:oscillatory-convex} with several differences
that we will point out here. Again, by interpolation,
it is sufficient to prove estimate
\[
 \Vert T_t f \Vert_{L^\infty(\Rn)} \leq 
 C \vert t\vert^{-\frac{n-1}{\gamma}}
 \Vert f \Vert_{{L}^1(\Rn)},
\]
with amplitude $a(t,\xi)$ satisfying
\begin{equation} \label{EQ:assAMP-low-freq1}
|\partial_\xi^\alpha a(t,\xi)|\leq C_\alpha |\xi|^{-N_1-|\alpha|}
\quad\text{for all}
\quad |\alpha|\leq [(n-1)/\gamma]+1,
\end{equation}
with $N_1=n-\frac{n-1}{\gamma}
+\left[ \frac{n-1}{\gamma} \right]+1$.

Further differences concern estimates for $I_1(t,x)$ and $I_2(t,x)$.
In general, since we work with low frequencies $|\xi|<1$ only,
no Besov space decomposition is necessary, so we do not need
to introduce function $\Psi$ and $\Phi_j$, so we can take
$\Psi=1$.

Some additional complications are related to the fact that 
in principle derivatives of the amplitudes of operators
$T_t$ from \eqref{EQ:oscillatory-cor} may introduce
an additional growth with respect to $t$.
In the estimate for $I_2(t,x)$ we performed integration
by parts with operator $P$. Now after integration by parts
the amplitude of this integral in \eqref{EQ:I2tx} is
$$
(P^*)^l\br{\br{1-\psi((1+|t|)\xi)} a(t,\xi) (1-\kappa)
\left(t^{-1}x+t^{-1}\nabla_\xi\vartheta(t;\xi)\right)
\Psi\p{\frac{\widetilde\vartheta(t;\xi)}{2^j}}}.
$$
Now, if any of the $\xi$-derivatives falls on 
$\br{1-\psi((1+|t|)\xi)}$, we get an extra factor $t$ which
is cancelled with $t^{-1}$ in the definition of $P$.
However, in this case we can then restrict to the support of
$\nabla\psi$ which is contained in the ball with radius
$(1+|t|)^{-1}$, so we are in the situation of
low frequencies $|\xi|\leq t^{-1}$ again. 
Consequently, we can apply Proposition
\ref{prop:low-freq} to this integral to actually get a 
better decay rate of Proposition \ref{prop:low-freq}.
If none of the derivatives in $(P^*)^l$ fall on
$\br{1-\psi((1+|t|)\xi)}$, the argument is the same as
in the proof of the estimate for $I_2(t,x)$ in
Proposition \ref{PROP:oscillatory-convex}.

The other main difference with the proof of Proposition
\ref{PROP:oscillatory-convex} is in the estimate
for $I_1(t,x)$. Recall now that in formula \eqref{EQ:I1tx}
we made a change of variables $\widetilde\lambda=x_n^{-1}\lambda$.
As it was then pointed out, if $r$ in the definition of
the cut-off function $\kappa$ is chosen sufficiently small,
on its support we have $|x_n|\approx |t|$. On the
other hand, we have $|\xi|\approx \widetilde\lambda$ 
by the definition
of $\widetilde\lambda$, since we assume that the limiting 
phase function $\va$ is strictly positive. It then follows
that $(1+|t|)\xi\approx \widetilde\lambda |x_n|\approx\lambda$, 
and so the change of variables $\widetilde\lambda=x_n^{-1}\lambda$
changes $\br{1-\psi((1+|t|)\xi)}$ into 
$\br{1-\psi(\lambda)}$ in the amplitude of $I_1(t,x)$.
Justifying this argument, we can then continue as in
the proof of Proposition \ref{PROP:oscillatory-convex}.
The crucial condition for the use of Theorem
\ref{THM:oscintthm} is the boundedness of derivatives of
$\widetilde{a}$ in \eqref{eq:a-tilde}. Here, every
differentiation of $a$ with respect to $y$ introduces
a factor $x_n^{-1}\lambda$ which is then cancelled in view of
assumption \eqref{EQ:assAMP-low-freq1}. It follows that
$[(n-1)/\gamma]+1$ $y$-derivatives of $\widetilde{a}$
are bounded, implying the conclusion of
Theorem \ref{THM:oscintthm}. This yields 
estimate \eqref{EQ:convex-estimate-cor} in the way that is
similar to the
proof of Proposition \ref{PROP:oscillatory-convex}.
\end{proof}

Let us now turn to prove Proposition 
\ref{PROP:oscillatory-nonconvex}.
\begin{proof}
[Proof of Proposition {\rm \ref{PROP:oscillatory-nonconvex}}]
Let us show how the proof of Proposition
\ref{PROP:oscillatory-nonconvex} differs from the proof
of Proposition \ref{PROP:oscillatory-convex}. 
We need to prove that $|I(t,x)|\leq Ct^{-\frac{1}{\gamma_0}}$, 
$t>0$, for $I(t,x)$ as in \eqref{EQ:Itx}. We note that 
$\gamma_0\geq 1$, so to prove the estimate for $I_2(t,x)$
we can show that $|I_2(t,x)|\leq Ct^{-1}.$ 
This can be done by integrating by parts with the same
operator $P$ and using \eqref{EQ:assAMP-nc} instead of
\eqref{EQ:assAMP}. As for the proof of the estimate
for $I_1(t,x)$, we can reason in the same way as in
Proposition \ref{PROP:oscillatory-convex} to arrive at the
estimate \eqref{EQ:IJ}, i.e.,
\[
|I_1(t,x)|\leq 
C t^{-n+l}2^{j(n-l)}\int_0^\infty 
\abs{J(2^j \lambda,z_t)
\Psi\p{\frac{\lambda}{t }}
\lambda^{n-1-l}} \, d\lambda
\]
with the same operator
$$
J(\lambda,z_t)=\int_U 
e^{i\lambda(-\nabla_y h_t(z_t)\cdot y+ h_t(y)+t x_n^{-1})}
\widetilde{a}(t,x_n,\lambda y,\lambda h_t(y)) 
\kappa_0(t,x,y) \chi(t,y) \, dy.
$$
Now, instead of \eqref{EQ:estJ} we will show that
\begin{equation}\label{EQ:estJ-nc}
\abs{J(\lambda,z_t)}\leq C\lambda^{-\frac{1}{\gamma_0}},
\quad \lambda>0.
\end{equation}
Then, taking $l=n-\frac{1}{\gamma_0}$, we get the estimate
$|I_1(t,x)|\leq C t^{-\frac{1}{\gamma_0}}$ in the same way
as in estimate \eqref{EQ:I1tx-est2}. Now, estimate
\eqref{EQ:estJ-nc} follows from Theorem \ref{THM:oscintthm}
with $N=1$. Indeed, let us write $J(\lambda,z_t)$
in the form \eqref{EQ:gather1}--\eqref{EQ:gather3} with
phase
$$F(\rho,z_t,\omega)= h_t(\rho\omega+z_t)-h_t(z_t)-
\rho \nabla_z h_t(z_t)\cdot\omega.$$ 
Now, by rotation we may assume that in some direction,
say $e_1=(1,0,\ldots,0)$, we have
by definition of the index $\gamma_0$ that
$$
\gamma_0=\min\brfig{k\in\N: \partial_{\omega_1}^k F(\rho,z_t,\omega)|_
{\omega_1=0}\not=0}.
$$
Then by taking $N=1$ and $y=\omega_1$ in Theorem 
\ref{THM:oscintthm}, we get the required estimate
\eqref{EQ:estJ-nc}.
\end{proof}

Now we will state the corollary of the proof of Proposition
\ref{PROP:oscillatory-nonconvex} which is similar to
Proposition \ref{COR:oscillatory-convex} to ensure its
application to our Cauchy problem. The proof is similar
to the proof of Proposition \ref{COR:oscillatory-convex}.
\begin{prop} \label{COR:oscillatory-nonconvex}
Let $T_t$, $t\ne0$, be an operator defined by
\[
T_t f(x)=\int_\Rn e^{i(x\cdot\xi+\vartheta(t;\xi))}
\br{1-\psi((1+|t|)\xi)}
a(t,\xi) \widehat{f}(\xi) \, d\xi,
\]
where $\vartheta(t;\xi)$ and 
$\gamma_0$ are as in Proposition
\ref{PROP:oscillatory-nonconvex}.
Assume that the amplitude
$a(t,\xi)$  satisfies $a(t,\xi)=0$ for all $|\xi|\geq 1$ and 
that 
\[
|\partial_\xi^\alpha a(t,\xi)|\leq C_\alpha |\xi|^{-|\alpha|}
\quad\text{for all}
\quad |\alpha|\leq 1.
\]
Let $1< p\leq 2\leq q< \infty$ be such that
$\frac{1}{p}+\frac{1}{q}=1$.
Then for $t\not=0$ we have the estimate
\begin{equation}\label{EQ:convex-estimate-cor-nc}
 \Vert T_t f \Vert_{L^q(\Rn)} \leq 
 C \vert t \vert^{-\frac{1}{\gamma_0}\p{\frac{1}{p}-
 \frac{1}{q}}}
 \Vert f \Vert_{\dot{L}^p_{N_p}(\Rn)},
\end{equation}
where $N_p=\left(n-\frac{1}{\gamma_0}+1\right)
\left(\frac{1}{p}-\frac{1}{q}\right)$. 
\end{prop}

Let us finally estimate our Fourier multiplier for small $t$. 
For the analysis of very small $t$ 
we will use the following Littlewood-Paley type theorem. 

\begin{lem}[\cite{Hormander}(Theorem~1.11)]
\label{lem:Hardy} 
Let $h=h(\xi)$ be a tempered distribution on $\mathbb{R}^{n}$ $(n \ge 1)$ 
such that 
\[
\sup_{0<l<+\infty} 
l^b {\rm meas}\, \{ {} \xi {} : {} \vert h(\xi) \vert \ge l {} \} 
<+\infty,
\]
for some $1<b<+\infty$. Then the convolution operator with 
$\mathcal{F}^{-1}[h]$ is $L^p$--$L^q$ bounded provided that 
$1<p \le 2 \le q<+\infty$ and 
$\frac{1}{p}-\frac{1}{q}=\frac{1}{b}$, i.e. we have the estimate
\[
\left\Vert \mathcal{F}^{-1}[h] *u \right\Vert_{L^q} 
\le C \Vert u \Vert_{L^p}\quad\textrm{for all} 
\quad u \in L^p(\mathbb{R}^n). 
\]
\end{lem}

Using this fact, we obtain: 

\begin{prop} \label{PROP:small time} 
Let $T_t$ be an operator defined by
\[
T_t f(x)=\int_\Rn e^{i(x\cdot\xi+\vartheta(t;\xi))}
\br{1-\psi((1+|t|)\xi)}
a(t,\xi) \widehat{f}(\xi) \, d\xi,
\]
where $\vartheta(t;\xi)$ and $a(t,\xi)$ 
are as in Propositions \ref{PROP:oscillatory-convex} or 
\ref{PROP:oscillatory-nonconvex}. Let $n \ge 1$, 
$1<p\le 2 \le q<+\infty$ and 
$\frac{1}{p}+\frac{1}{q}=1$. 
Then for small $t$ we have the estimate 
\begin{equation}
\Vert T_t f \Vert_{L^q(\mathbb{R}^n)} \le C 
\Vert f \Vert_{\dot{L}^{p}_{\widetilde{N}_p}(\mathbb{R}^n)}, 
\label{another estimates} 
\end{equation} 
where $\widetilde{N}_p=n \left(\frac{1}{p}-\frac{1}{q} \right)$.
\end{prop} 
In fact, the proof will yield the Besov norm $B_p^{1,1}$ on the
right hand side of the estimate, which is a known improvement
for this type of estimates.
\begin{proof}
In the following argument we need not use 
the Van der Corpt lemma, 
and the proof relies only on the Littlewood--Paley type 
theorem. 
We put $K(t)=(1+|t|)^{-1}$. 

It suffices to prove \eqref{another estimates} for $p\ne q$, 
since the case $p=q=2$ follows from the 
Plancherel theorem. Noting that $\vartheta(t;\xi)$ and is 
homogeneous of order one, and making change of variable
$\eta=\frac{\xi}{K(t)}$ and $y=K(t)x$, we get  
\begin{equation}
\Vert T_t f \Vert_{L^q(\mathbb{R}^n)}= 
K(t)^{n-\frac{n}{q}}\left\Vert 
\mathcal{F}^{-1} \left[m_t(\eta)\right]
* \mathcal{F}^{-1} \left[\vert \eta \vert^{\widetilde{N}_p}
\widehat{f}(K(t)\eta)\right] 
\right\Vert_{L^q(\mathbb{R}^n)},
\label{Paley}
\end{equation}
where we set 
\[
m_t(\eta)=e^{iK(t)\vartheta(t;\eta)}a(t,K(t)\eta)
(1-\psi(\eta))|\eta|^{-\widetilde{N}_p}, \quad 
\widetilde{N}_p=n \left(\frac{1}{p}-\frac{1}{q} \right).
\]
Since $a(t,\xi)$ is bounded, we have 
\[
{\rm meas} {} \left\{ {} \eta {} : {} |m_t(\eta)| 
\ge l {} \right\} \le
{\rm meas} {} \left\{ {} \eta {} : {} |\eta| \le 
C^{1/\widetilde{N}_p}l^{-1/\widetilde{N}_p} 
{} \right\}
=C^{n/\widetilde{N}_p}l^{-n/\widetilde{N}_p} 
\]
for each $l>0$. Hence it follows from Lemma 
\ref{lem:Hardy} that 
the convolution operator with $m_t$ is $L^p$--$L^q$ bounded, 
which implies that 
\begin{eqnarray*}
\left\Vert \mathcal{F}^{-1}[m_t] * \mathcal{F}^{-1} 
\left[\vert \eta \vert^{\widetilde{N}_p}\widehat{f}(K(t)\eta)\right] 
\right\Vert_{L^q(\mathbb{R}^n)} 
&\le& C \left\Vert \mathcal{F}^{-1} \left[ \vert \eta \vert^{\widetilde{N}_p}
\widehat{f}(K(t)\eta)\right] \right\Vert_{L^p(\mathbb{R}^n)}\\ 
&=&C K(t)^{-n+\frac{n}{p}-\widetilde{N}_p} 
\Vert f \Vert_{\dot{L}^{p}_{\widetilde{N}_p}(\mathbb{R}^n)},
\end{eqnarray*}
where we performed the transformations $K(t)\eta=\xi$ and 
$\frac{x}{K(t)}=z$ in the last step. Thus, combining this 
estimate with \eqref{Paley}, we obtain the desired 
estimate \eqref{another estimates}. The proof of 
Proposition~\ref{PROP:small time} is complete. 
\end{proof} 
\begin{proof}
[Proof of Theorem {\rm \ref{thm:Main theorem}}]
The proof of Theorem \ref{thm:Main theorem} now follows
from Proposition \ref{prop:low-freq} for low frequencies
$|\xi|<t^{-1}$, from Propositions
\ref{PROP:oscillatory-convex} and
\ref{PROP:oscillatory-nonconvex} for large frequencies
$|\xi|\geq 1$, and from Propositions
\ref{COR:oscillatory-convex} and \ref{COR:oscillatory-nonconvex}
for intermediate frequencies $t^{-1}\leq |\xi|<1$.
We also use Proposition \ref{PROP:small time} for small times.
We can note that all these propositions give different Sobolev
orders on the regularity of the Cauchy data.

Indeed, using representation formula for the solution established
in Theorem \ref{thm:asymptotic}, we can write the solution
as $$u(t,x)=\sum_{k=0}^{m-1} u^{k}(t,x),$$ 
with
$$
  u^{k}(t,x)=\sum_{j=1}^m 
  \mathcal{F}^{-1} 
\left[ \left(\alpha^{j}_{k,\pm}(\xi)
+\varepsilon^{j}_{k,\pm}(t;\xi)\right)  
e^{i\vartheta_j(t;\xi)} \widehat{f}_k (\xi)\right](x).
$$
Now, we decompose
\begin{align*}
& u^{k}(t,x)= u^{k}_1(t,x) + u^{k}_2(t,x)+u^{k}_3(t,x) \\
 & =
 \sum_{j=1}^m 
  \mathcal{F}^{-1} 
\left[ \left(\alpha^{j}_{k,\pm}(\xi)
+\varepsilon^{j}_{k,\pm}(t;\xi)\right) \psi((1+|t|)\xi) 
e^{i\vartheta_j(t;\xi)} \widehat{f}_k (\xi)\right](x)\\
& + 
\sum_{j=1}^m 
  \mathcal{F}^{-1} 
\left[ \left(\alpha^{j}_{k,\pm}(\xi)
+\varepsilon^{j}_{k,\pm}(t;\xi)\right) (1-\psi((1+|t|)\xi))\chi(\xi)
e^{i\vartheta_j(t;\xi)} \widehat{f}_k (\xi)\right](x)\\
& + \sum_{j=1}^m 
  \mathcal{F}^{-1} 
\left[ \left(\alpha^{j}_{k,\pm}(\xi)
+\varepsilon^{j}_{k,\pm}(t;\xi)\right) (1-\psi((1+|t|)\xi))
(1-\chi(\xi)) 
e^{i\vartheta_j(t;\xi)} \widehat{f}_k (\xi)\right](x),
\end{align*}
with
$\psi,\chi\in C_0^\infty(\Rn)$
such that $\psi(\xi)=\chi(\xi)\equiv 1$ for 
$\vert \xi \vert \le \frac12$, 
and $0$ for  $\vert \xi \vert \ge 1$. 
Assume conditions of part (i) of Theorem
\ref{thm:Main theorem}. Then we have estimates
\begin{equation}
  \Vert u^k_1(t,\cdot)\Vert_{L^q(\Rn)} \leq 
 C(1+\vert t\vert)^{-n\left(\frac{1}{p}-\frac{1}{q}\right)}
 \Vert f_k \Vert_{\dot{L}^{p}_{-k}(\Rn)}, 
\label{solest1} 
\end{equation}
by Proposition \ref{prop:low-freq},
\begin{equation}
 \Vert u^k_2(t,\cdot)  \Vert_{L^q(\Rn)} \leq 
 C \vert t\vert^{-\frac{n-1}{\gamma}\p{\frac{1}{p}-
 \frac{1}{q}}}
 \Vert f_k \Vert_{\dot{L}^p_{M_p-k}(\Rn)}, 
\label{solest2} 
\end{equation}
with $M_p=\left(n-\frac{n-1}{\gamma}
\right)
\left(\frac{1}{p}-\frac{1}{q}\right)$
by Proposition \ref{COR:oscillatory-convex}, and
\begin{equation}
\Vert u^k_3(t,\cdot) \Vert_{L^q(\Rn)} 
 \leq C |t|^{-\frac{n-1}{\gamma}\p{\frac{1}{p}-\frac{1}{q}}}
 \Vert f_k \Vert_{{L}^p_{N_p-k}(\Rn)},
\label{solest3} 
\end{equation}
with $N_p=\left(n-\frac{n-1}{\gamma}+
\left[ \frac{n-1}{\gamma} \right]+1\right)
\left(\frac{1}{p}-\frac{1}{q}\right)$,
by Proposition \ref{PROP:oscillatory-convex}.
For small $t$ we have the estimate
\begin{equation}
\Vert u^k_2(t,\cdot) \Vert_{L^q(\mathbb{R}^n)} 
+ \Vert u^k_3(t,\cdot) \Vert_{L^q(\mathbb{R}^n)}\le C 
\Vert f_k \Vert_{\dot{L}^{p}_{\widetilde{N}_p-k}(\mathbb{R}^n)}, 
\label{solest4} 
\end{equation} 
with $\widetilde{N}_p=n \left(\frac{1}{p}-\frac{1}{q} \right)$,
by Proposition \ref{PROP:small time}.
Putting all these estimates \eqref{solest1}--\eqref{solest4} 
together with similar estimates 
for derivatives, implies the statement of Theorem 
\ref{thm:Main theorem}.
\end{proof}
As a corollary of this proof and all the propositions above,
we have the following refinement of Theorem \ref{thm:Main theorem},
providing quantitatively different estimates 
for different frequency regions.

\begin{thm}\label{thm:Main theorem-2}
Assume \eqref{L1}--\eqref{strict hyperbolicity2}. 
Let $\chi\in C_0^\infty(\R)$ be
such that $\chi(\rho)\equiv 1$ for 
$\vert \rho \vert \le \frac12$, 
and $0$ for  $\vert \rho \vert \ge 1$. 
Let us denote
$$u_1=\chi((1+|t|)|D|)u, \quad 
u_2=(1-\chi((1+|t|)|D|)\chi(|D|)u, \quad \mathrm{and}$$
$$u_3=(1-\chi((1+|t|)|D|)(1-\chi(|D|))u.$$
Then the solution $u(t,x)$ of \eqref{Equation} 
satisfies the following 
estimates{\rm :}

{\rm (i)} Suppose that the sets
$$\Si_{\varphi^{\pm}_k}=
\{\xi\in\Rn: \va^{\pm}_k(\xi)=1\}$$ are
convex for all $k=1,\ldots,m$, and set
$\gamma=\underset{k=1,\ldots,m}{\max} 
\gamma(\Sigma_{\varphi^{\pm}_k})$. 
In addition, suppose that 
$(1+\vert t \vert)^r a^\prime_{\nu,j} 
\in L^1(\mathbb{R})$ for 
$1\le r \le [(n-1)/\gamma]+1$, 
and for all $\nu,j$ with $|\nu|+j=m$. 
Let $1<p\le 2 \le q<+\infty$ and $\frac{1}{p}+\frac{1}{q}=1$. 
Then we have the estimates
\begin{gather*}
\Vert D_t^l D_x^\alpha u_1(t,\cdot)\Vert_{L^q(\Rn)} \leq 
 C(1+\vert t\vert)^{-n\left(\frac{1}{p}-\frac{1}{q}\right)}
 \sum_{k=0}^{m-1} 
 \Vert f_k \Vert_{\dot{L}^{p}_{l+|\alpha|-k}(\Rn)}, 
\quad (t\in \mathbb{R}),\\
 \Vert D_t^l D_x^\alpha u_2(t,\cdot)  \Vert_{L^q(\Rn)} \leq 
 C \vert t\vert^{-\frac{n-1}{\gamma}\p{\frac{1}{p}-
 \frac{1}{q}}} \sum_{k=0}^{m-1}
 \Vert f_k \Vert_{\dot{L}^p_{l+|\alpha|+M_p-k}(\Rn)}, \quad 
 (|t|\ge1),
\\
\Vert D_t^l D_x^\alpha u_3(t,\cdot) \Vert_{L^q(\Rn)} 
 \leq C |t|^{-\frac{n-1}{\gamma}\p{\frac{1}{p}-\frac{1}{q}}}
 \sum_{k=0}^{m-1}
 \Vert f_k \Vert_{{L}^p_{l+|\alpha|+N_p-k}(\Rn)}, \quad (|t|\ge1),\\
\Vert D_t^l D_x^\alpha u_2(t,\cdot) \Vert_{L^q(\mathbb{R}^n)} 
+ \Vert D_t^l D_x^\alpha 
u_3(t,\cdot) \Vert_{L^q(\mathbb{R}^n)} \le C 
\sum_{k=0}^{m-1}
\Vert f_k \Vert_{\dot{L}^{p}_{l+|\alpha|+
\widetilde{N}_p-k}(\mathbb{R}^n)}, \quad 
(|t|<1),
\end{gather*}
with $M_p=\left(n-\frac{n-1}{\gamma}\right)
\left(\frac{1}{p}-\frac{1}{q}\right)$, 
$N_p=\left(n-\frac{n-1}{\gamma}+
\left[ \frac{n-1}{\gamma} \right]+1\right)
\left(\frac{1}{p}-\frac{1}{q}\right)$,
$\widetilde{N}_p=n \left(\frac{1}{p}-\frac{1}{q} \right)$,
$l=0,\ldots,m-1$, and 
$\alpha$ any multi-index.

{\rm (ii)} Suppose that $\Si_{\varphi^{\pm}_k}$ 
is non-convex 
for some $k=1,\ldots,m$, and set $\gamma_0=
\underset{k=1,\ldots,m}{\max} 
\gamma_0(\Sigma_{\varphi^{\pm}_k})$. 
In addition, suppose that 
$(1+\vert t \vert) a^\prime_{\nu,j} \in L^1(\mathbb{R})$ 
for all $\nu,j$ with $|\nu|+j=m$. Let $1<p\le 2 \le q<+\infty$ 
and $\frac{1}{p}+\frac{1}{q}=1$. 
Then we have the estimates 
\begin{gather*}
  \Vert D_t^l D_x^\alpha u_1(t,\cdot)\Vert_{L^q(\Rn)} \leq 
 C(1+\vert t\vert)^{-n\left(\frac{1}{p}-\frac{1}{q}\right)}
 \sum_{k=0}^{m-1} 
 \Vert f_k \Vert_{\dot{L}^{p}_{l+|\alpha|-k}(\Rn)}, \quad (t\in \mathbb{R}),
\\
 \Vert D_t^l D_x^\alpha u_2(t,\cdot)  \Vert_{L^q(\Rn)} \leq 
 C \vert t\vert^{-\frac{1}{\gamma_0}\p{\frac{1}{p}-
 \frac{1}{q}}} \sum_{k=0}^{m-1}
 \Vert f_k \Vert_{\dot{L}^p_{l+|\alpha|+M_p-k}(\Rn)}, \quad (|t|\ge1),
\\
 \Vert D_t^l D_x^\alpha u_3(t,\cdot) \Vert_{L^q(\Rn)} 
 \leq C |t|^{-\frac{1}{\gamma_0}\p{\frac{1}{p}-\frac{1}{q}}}
 \sum_{k=0}^{m-1}
 \Vert f_k \Vert_{{L}^p_{l+|\alpha|+N_p-k}(\Rn)}, \quad (|t|\ge1),
\\
\Vert D_t^l D_x^\alpha u_2(t,\cdot) \Vert_{L^q(\mathbb{R}^n)} 
+ \Vert D_t^l D_x^\alpha 
u_3(t,\cdot) \Vert_{L^q(\mathbb{R}^n)} \le C 
\sum_{k=0}^{m-1}
\Vert f_k \Vert_{\dot{L}^{p}_{l+|\alpha|+
\widetilde{N}_p-k}(\mathbb{R}^n)}, \quad 
(|t|<1),
\end{gather*}
with $M_p=\left(n-\frac{1}{\gamma_0}\right)
\left(\frac{1}{p}-\frac{1}{q}\right)$, 
$N_p=\left(n-\frac{1}{\gamma_0}+1\right)
\left(\frac{1}{p}-\frac{1}{q}\right)$,
$\widetilde{N}_p=n \left(\frac{1}{p}-\frac{1}{q} \right)$,
$l=0,\ldots,m-1$, and 
$\alpha$ any multi-index.
\end{thm}



\begin{thebibliography}{99} 

\bibitem{AVG} V.~I.~Arnold, S.~M.~Gusein-Zade and
A.~N.~Varchenko, ``Singularities of differentiable maps''. 
Vol. I. 
The classification of critical points, caustics and 
wave fronts. Monographs in Mathematics, 82. 
Birkh\"auser Boston, Inc., Boston, MA, 1985.



%
\bibitem{Ascoli}G.~Ascoli, 
{\em Sulla forma asintotica degli integrali 
dell'equazione differenziale 
$y^{\prime \prime}+A(x)y=0$ in un 
caso notevole di stabilit\`a}, Univ. Nac. 
Tucum\'an, Revista A. {\bf 2} (1941), 131--140.

\bibitem{BL} J.~Bergh and J.~L\"ofstr\"om,
 ``Interpolation spaces'', Springer, 1976.


\bibitem{Bren75}
P.~Brenner,  
{\em On {$L\sb{p}$--$L\sb{p\sp{\prime} }$} estimates for the
  wave equation}, Math. Z. \textbf{145} (1975), 251--254.

\bibitem{Bren77}
P.~Brenner, 
{\em $L\sb{p}$--$L\sb{p'}$-estimates for {F}ourier integral operators
  related to hyperbolic equations}, Math. Z. 
  \textbf{152} (1977), 273--286.


%
\bibitem{Coddington}
   E.A. Coddington and N. Levinson, 
``Theory of differential equations'', 
New York, McGraw-Hill, 1955. 


%
 \bibitem{Hartman} 
  P.~Hartman, ``Ordinary differential equations'', 
  SIAM, 2nd edition, 2002. 
%

\bibitem{HR}
F.~Hirosawa and M.~Reissig, {\em 
About the optimality of oscillations in non-Lipschitz 
coefficients for  strictly hyperbolic equations}, 
Ann. Sc. Norm. Super. Pisa Cl. Sci. (5) {\bf 3} (2004), 
589--608.
                
%
\bibitem{Hormander}L.~H\"ormander, 
{\em  Estimates for translation 
invariant operators in $L\sp{p}$ spaces}, 
Acta Math. {\bf 104} (1960), 93--140.



\bibitem{Litt73}
W.~Littman, 
{\em {$L\sp{p}$--$L\sp{q}$}-estimates for singular integral operators
  arising from hyperbolic equations},
   Partial differential equations (Proc.
  Sympos. Pure Math., Vol. XXIII, Univ. California, 
  Berkeley, Calif., 1971),
  Amer. Math. Soc., Providence, R.I., 1973, pp.~479--481.


\bibitem{Matsuyama0}T.~Matsuyama, 
{\em Asymptotic behaviour for wave equations 
with time-dependent coefficients}, Annali dell'Universit\`a di Ferrara, 
Sec. VII  - Sci. Math., {\bf 52} (2), (2006), 383--393.

\bibitem{Matsuyama1}T.~Matsuyama, 
{\em $L^p$--$L^q$ estimates for wave 
equations and the Kirchhoff equation}, 
Osaka J. Math. {\bf 45} (2008), 491--510. 

\bibitem{Matsuyama2} T.~Matsuyama and M.~Reissig, 
{\em Stabilization and 
$L^p$--$L^q$ decay estimates}, Asymptotic Anal. \textbf{50} (2006), 239--268.

\bibitem{MR09}
T.~Matsuyama and M.~Ruzhansky M., 
{\em Time decay for hyperbolic equations with homogeneous 
symbols}, C. R. Acad. Sci. Paris, Ser  I. {\bf 347} (2009), 915--919.
 
%
\bibitem{Mizohata}S.~Mizohata, 
``The theory of partial differential 
equations'', Cambridge Univ. Press, 1973. 

\bibitem{Reissig1} M.~Reissig, {\em $L\sb p$--$L\sb q$ decay 
estimates for wave equations with time-dependent coefficients}, 
J. Nonlinear Math. Phys. {\bf 11} (2004), 534--548.
%
\bibitem{Reissig} 
M.~Reissig and J.~Smith, {\em $L^p$--$L^q$ estimates 
for wave equation with 
bounded time dependent coefficient}, Hokkaido Math. J. 
{\bf 34} (2005), 541--586. 

\bibitem{Ruzh-survey} M.~Ruzhansky,
{\em Singularities of affine fibrations in the theory 
of regularity of Fourier integral operators}, 
Russian Math. Surveys {\bf 55} (2000), 93--161.

\bibitem{R} M.~Ruzhansky,
{\em On some properties of Galerkin approximations of solutions
to Fokker--Planck equations}, 
in Proceedings of the 4th  International Conference 
"Analytical Methods in Analysis and Differential Equations" 
(AMADE-2006), Vol.3, Differential Equations, Minsk: 
Institute of Mathematics of NAS of Belarus, 133-139, 2006.

\bibitem{R-vdC} M.~Ruzhansky,
{\em Pointwise van der Corput lemma for functions 
of several variables}, Funct. Anal. and Appl.
{\bf 43} (2009), 75--77.

\bibitem{RS} M.~Ruzhansky and J.~Smith, 
{\em Global time estimates for solutions to equations of 
dissipative types,} 
Journees ``Equations aux Derivees Partielles", Exp. No. XII, 
29 pp., Ecole Polytech., Palaiseau, 2005.

\bibitem{RS2} M.~Ruzhansky and J.~Smith, 
``Dispersive and Strichartz estimates for 
hyperbolic equations with constant coefficients'', 
arXiv:0711.2138v1, to appear in MSJ Memoirs, Vol.22, 2010. 

\bibitem{stri70func}
R.~Strichartz, 
\emph{A priori estimates for the wave equation and some applications},
  J. Func. Anal. \textbf{5} (1970), 218--235.

\bibitem{sugi94}
M.~Sugimoto, 
\emph{A priori estimates for higher order hyperbolic equations},
  Math. Z. \textbf{215} (1994), 519--531.

\bibitem{sugi96}
M.~Sugimoto, \emph{Estimates for hyperbolic equations with non-convex
  characteristics}, Math. Z. \textbf{222} (1996), 521--531.

\bibitem{sugi98}
M.~Sugimoto, 
\emph{Estimates for hyperbolic equations of space dimension 3}, J.
  Funct. Anal. \textbf{160} (1998), 382--407.

%
 \bibitem{Wintner} 
  A.~Wintner, 
  {\em Asymptotic integrations of adiabatic oscillator}, 
Amer.\ J.\ Math.\ {\bf 69} (1947), 251--272. 

%
 
%


\end{thebibliography}
\end{document}